\documentclass[12pt]{article}%
\usepackage{amsmath}
\usepackage{amsfonts}
\usepackage{amssymb}
\usepackage{graphicx}%
\providecommand{\U}[1]{\protect\rule{.1in}{.1in}}
\begin{document}

\title{Projections over Quantum Homogeneous Odd-dimensional Spheres\thanks{This work
was partially supported by the grant H2020-MSCA-RISE-2015-691246-QUANTUM
DYNAMICS and the Polish government grant 3542/H2020/2016/2.}}
\date{}
\author{Albert Jeu-Liang Sheu\\{\small Department of Mathematics, University of Kansas, Lawrence, KS 66045,
U. S. A.}\\{\small e-mail: asheu@ku.edu}}
\maketitle

\begin{abstract}
We give a complete classification of isomorphism classes of finitely generated
projective modules, or equivalently, unitary equivalence classes of
projections, over the C*-algebra $C\left(  \mathbb{S}_{q}^{2n+1}\right)  $ of
the quantum homogeneous sphere $\mathbb{S}_{q}^{2n+1}$. Then we explicitly
identify as concrete elementary projections the quantum line bundles $L_{k}$
over the quantum complex projective space $\mathbb{C}P_{q}^{n}$ associated
with the quantum Hopf principal $U\left(  1\right)  $-bundle $\mathbb{S}%
_{q}^{2n+1}\rightarrow\mathbb{C}P_{q}^{n}$.

\end{abstract}

\pagebreak

\section{Introduction}

In the theory of quantum/noncommutative geometry popularized by Connes
\cite{Conn}, C*-algebras are often viewed as the algebra $C\left(
X_{q}\right)  $ of continuous functions on a virtual quantum space $X_{q}$,
and finitely generated projective (left) $C\left(  X_{q}\right)  $-module
$\Gamma\left(  E_{q}\right)  $ are viewed as virtual vector bundles over the
quantum space $X_{q}$. The former viewpoint is motivated by Gelfand's Theorem
identifying all commutative C*-algebras as exactly function algebras
$C_{0}\left(  X\right)  $ for locally compact Hausdorff spaces $X$, while the
latter is motivated by Swan's Theorem \cite{Swan} characterizing all finitely
generated projective $C\left(  X\right)  $-modules for a compact Hausdorff
space $X$ as exactly the spaces $\Gamma\left(  E\right)  $ of continuous
cross-sections of vector bundles $E$ over $X$.

As spheres and projective spaces provide fundamentally important examples for
the classical study of topology and geometry, quantum versions of spheres and
projective spaces have been developed and provide important examples for the
study of quantum geometry. In particular, from the quantum group viewpoint
\cite{Dr,Woro,Wo:cm}, Soibelman, Vaksman, Meyer and others
\cite{VaSo,Me,Na,Sh:qcp} introduced and studied the homogeneous
odd-dimensional quantum sphere $\mathbb{S}_{q}^{2n+1}$ and the associated
quantum complex projective space $\mathbb{C}P_{q}^{n}$, and from the
multipullback viewpoint, Hajac and his collaborators including Baum, Kaygun,
Matthes, Nest, Pask, Sims, Szyma\'{n}ski, Zieli\'{n}ski, and others
\cite{HaMaSz,HaKaZi,HaNePaSiZi} developed and studied the multipullback
odd-dimensional quantum sphere $\mathbb{S}_{H}^{2n+1}$ and the associated
quantum complex projective space $\mathbb{P}^{n}\left(  \mathcal{T}\right)  $.

As in the classical situation, the above mentioned quantum odd-dimensional
spheres and their associated quantum complex projective spaces provide a
quantum Hopf principal $U\left(  1\right)  $-bundle, from which some
associated quantum line bundles $L_{k}$, or rank-one projective modules over
the quantum algebra of the complex projective space, for $k\in\mathbb{Z}$ are
constructed \cite{Me,AriBrLa,HaMaSz,HaNePaSiZi}.

It is well known that classifying up to isomorphism all vector bundles over a
space $X$ in the classical case or finitely generated (left) projective
modules over a C*-algebra $C\left(  X_{q}\right)  $ in the quantum case is an
interesting but difficulty task. A major challenge in such classification is
the so-called cancellation problem \cite{Ri:dsr,Ri:ct} which deals with
determining whether the stable isomorphism between such objects determined by
$K$-theoretic analysis can imply their isomorphism.

In this paper, we use the powerful groupoid approach to C*-algebras initiated
by Renault \cite{Rena} and popularized by Curto, Muhly, and Renault
\cite{CuMu,MuRe} to study the C*-algebra structures of $C\left(
\mathbb{S}_{q}^{2n+1}\right)  $ and $C\left(  \mathbb{C}P_{q}^{n}\right)  $.
In this framework, we work to get a complete classification of projections
over $C\left(  \mathbb{S}_{q}^{2n+1}\right)  $ up to equivalence, extending
the result of Bach \cite{Ba}, and determine the canonical monoid structure on
the collection of all equivalence classes of projections over $C\left(
\mathbb{S}_{q}^{2n+1}\right)  $ with the diagonal sum $\boxplus$ as its binary
operation. In particular, we get infinitely many inequivalent projections over
$C\left(  \mathbb{S}_{q}^{2n+1}\right)  $ which are stably equivalent over
$C\left(  \mathbb{S}_{q}^{2n+1}\right)  $, showing that the cancellation
property does not hold for projections over $C\left(  \mathbb{S}_{q}%
^{2n+1}\right)  $ as elaborated in Corollary 1. Then we proceed to present a
set of elementary projections that freely generate $K_{0}\left(  C\left(
\mathbb{C}P_{q}^{n}\right)  \right)  $ and represent the line bundles $L_{k}$
over $\mathbb{C}P_{q}^{n}$ by concrete $\boxplus$-sum of elementary
projections. We mention that a similar study has been carried out for the
multipullback quantum spheres $\mathbb{S}_{H}^{2n+1}$ and the associated
projective space $\mathbb{P}^{n}\left(  \mathcal{T}\right)  $ in the paper
\cite{Sh:qps,Sh:pmqpl}, and an interesting geometric study via Milnor
construction is presented by Farsi, Hajac, Maszczyk, and Zieli\'{n}ski in
\cite{FHMZ} for $C\left(  \mathbb{P}^{2}\left(  \mathcal{T}\right)  \right)  $.

Among works in the literature related to our topic here, we mention that the
graph C*-algebra of any row-finite graph, including $C\left(  \mathbb{S}%
_{q}^{2n+1}\right)  $, satisfies the so-called stable weak cancellation
property \cite{AraMoPa}, and that a \textquotedblleft
geometric\textquotedblright\ realization of generators of $K_{0}\left(
C\left(  \mathbb{C}P_{q}^{n}\right)  \right)  $ using Milnor connecting
homomorphism is found in \cite{ADHT}, beside the geometric study of quantum
line bundles over $\mathbb{C}P_{q}^{n}$ in \cite{AriBrLa}. It would be of
interest to take a close look at potential underlying connections between
these works and ours. (The author thanks the referee for relevant references
to the literature.)

The author would like to thank Prof. Dabrowski for hosting his visit to SISSA,
Trieste, Italy in the spring of 2018, and also thank him and Prof. Landi for
useful discussions and questions about quantum odd-dimensional spheres and
quantum complex projective spaces.

\section{Preliminary notations}

In this paper, we use freely the basic techniques and manipulations for
$K$-theory of C*-algebras, or more generally, Banach algebras, found in
\cite{Blac,Tay}. Commonly widely used notations like $M_{\infty}\left(
\mathcal{A}\right)  $, $GL_{\infty}\left(  \mathcal{A}\right)  $, unitization
$\mathcal{A}^{+}$, diagonal sum $P\boxplus Q$ of elements $P,Q\in M_{\infty
}\left(  \mathcal{A}\right)  $, the identity component $G^{0}$ of a
topological group $G$, the positive cone $K_{0}\left(  \mathcal{A}\right)
_{+}$ of $K_{0}\left(  \mathcal{A}\right)  $, $\mathcal{B}\left(
\mathcal{H}\right)  $, $\mathcal{K}\left(  \mathcal{H}\right)  $, and others
will not be explained in details here, and we refer to the notations section
in \cite{Sh:qps} for any need of further clarification.

By a projection (or an idempotent) over a C*-algebra $\mathcal{A}$, we mean a
projection (or an idempotent) in the algebra $M_{\infty}\left(  \mathcal{A}%
\right)  $ of all finite matrices with entries in $\mathcal{A}$. Two
projections (or idempotents) $P,Q\in M_{\infty}\left(  \mathcal{A}\right)  $
are called equivalent over $\mathcal{A}$, denoted as $P\sim_{\mathcal{A}}Q$,
if there is an invertible $U\in GL_{\infty}\left(  \mathcal{A}\right)  $ such
that $UPU^{-1}=Q$.

We recall that the mapping $P\mapsto\mathcal{A}^{n}P$ induces a bijective
correspondence between the equivalence (respectively, the stable equivalence)
classes of idempotents over $\mathcal{A}$ and the isomorphism (respectively,
the stable isomorphism) classes of finitely generated projective modules over
$\mathcal{A}$ \cite{Blac}, where by a module over $\mathcal{A}$, we mean a
left $\mathcal{A}$-module, unless otherwise specified.

We also recall that the $K_{0}$-group $K_{0}\left(  \mathcal{A}\right)  $
classifies idempotents over $\mathcal{A}$ up to stable equivalence. The
classification of idempotents up to equivalence, appearing as the so-called
cancellation problem, was popularized by Rieffel's pioneering work
\cite{Ri:dsr,Ri:ct} and is in general an interesting but difficult question.

For a C*-algebra homomorphism $h:\mathcal{A}\rightarrow\mathcal{B}$, we use
the same symbol $h$, instead of the more formal symbol $M_{\infty}\left(
h\right)  $, to denote the algebra homomorphism $M_{\infty}\left(
\mathcal{A}\right)  \rightarrow M_{\infty}\left(  \mathcal{B}\right)  $ that
applies $h$ to each entry of any matrix in $M_{\infty}\left(  \mathcal{A}%
\right)  $.

The set of all equivalence classes of idempotents, or equivalently, all
unitary equivalence classes of projections, over a C*-algebra $\mathcal{A}$ is
an abelian monoid $\mathfrak{P}\left(  \mathcal{A}\right)  $ with its binary
operation provided by the diagonal sum $\boxplus$.

In the following, we use the notations $\mathbb{Z}_{\geq k}:=\left\{
n\in\mathbb{Z}|n\geq k\right\}  $ and $\mathbb{Z}_{\geq}:=\mathbb{Z}_{\geq0}$.
In particular, $\mathbb{N}=\mathbb{Z}_{\geq1}$. We use $I$ to denote the
identity operator canonically contained in $\mathcal{K}^{+}\subset
\mathcal{B}\left(  \ell^{2}\left(  \mathbb{Z}_{\geq}\right)  \right)  $, and
\[
P_{m}:=\sum_{i=1}^{m}e_{ii}\in M_{m}\left(  \mathbb{C}\right)  \subset
\mathcal{K}%
\]
to denote the standard $m\times m$ identity matrix in $M_{m}\left(
\mathbb{C}\right)  \subset\mathcal{K}$ for any integer $m\geq0$ (with
$M_{0}\left(  \mathbb{C}\right)  =0$ and $P_{0}=0$ understood). We also use
the notation
\[
P_{-m}:=I-P_{m}\in\mathcal{K}^{+}%
\]
for integers $m>0$, and take symbolically $P_{-0}\equiv I-P_{0}=I\neq P_{0}$.
This should not cause any trouble since we will not formally add up the
subscripts of these $P$-projections without necessary clarification.

\section{Quantum spaces as groupoid C*-algebras}

In the following, we work with some concrete $r$-discrete (or \'{e}tale)
groupoids and use them to analyze and encode important structures of quantum
$\mathbb{S}_{q}^{2n+1}$ and quantum $\mathbb{C}P_{q}^{n}$ in the context of
groupoid C*-algebras. This groupoid approach to C*-algebras was popularized by
the work of Curto, Muhly, and Renault \cite{CuMu,MuRe,SSU} and shown to be
useful in the study of quantum homogeneous spaces
\cite{Sh:cqg,Sh:qpsu,Sh:qs,Sh:qcp}. We refer readers to Renault's pioneering
book \cite{Rena} for the fundamental theory of groupoid C*-algebras.

We denote by $\overline{\mathbb{Z}}:=\mathbb{Z}\cup\left\{  +\infty\right\}  $
the discrete space $\mathbb{Z}$ with a point $+\infty\equiv\infty$ canonically
adjoined as a limit point at the positive end, and take $\mathbb{Z}_{\geq
}:\equiv\left\{  n\in\mathbb{Z}|n\geq0\right\}  \subset\overline{\mathbb{Z}}$.
(We could also take $\overline{\mathbb{Z}}$ to be the one-point
compactification of the discrete space $\mathbb{Z}$ in this paper since
essentially we work only with groupoids restricted to a positive cone of their
unit spaces.) The group $\mathbb{Z}$ acts by homeomorphisms on $\overline
{\mathbb{Z}}$ in the canonical way, namely, by translations on $\mathbb{Z}$
while fixing the point $\infty$. More generally, the group $\mathbb{Z}^{n}$
acts on $\overline{\mathbb{Z}}^{n}$ componentwise in such a way. Let
$\mathcal{F}^{n}:={\mathbb{Z}}\times\left(  {\mathbb{Z}}^{n}\ltimes
\overline{{\mathbb{Z}}}^{n}\right)  |_{\overline{{\mathbb{Z}}}_{\geq}^{n}}$
with $n\geq1$ be the direct product of the group $\mathbb{Z}$ and the
transformation groupoid $\mathbb{Z}^{n}\ltimes\overline{\mathbb{Z}}^{n}$
restricted to the positive \textquotedblleft cone\textquotedblright%
\ $\overline{\mathbb{Z}_{\geq}}^{n}$, where $\overline{\mathbb{Z}_{\geq}}$ is
the closure $\mathbb{Z}_{\geq}\cup\left\{  \infty\right\}  $ of $\mathbb{Z}%
_{\geq}$ in $\overline{\mathbb{Z}}$. (Later we also use $\overline{\mathbb{Z}%
}_{\geq}$ to denote this positive part $\overline{\mathbb{Z}_{\geq}}$ of
$\overline{\mathbb{Z}}$.)

As shown in \cite{Sh:qs}, $C\left(  \mathbb{S}_{q}^{2n+1}\right)  \simeq
C^{\ast}\left(  \mathfrak{F}_{n}\right)  $, where $\mathfrak{F}_{n}$ is a
subquotient groupoid of $\mathcal{F}^{n}$, namely, $\mathfrak{F}%
_{n}:=\widetilde{\mathfrak{F}_{n}}/\sim$ for the subgroupoid%

\[
\widetilde{\mathfrak{F}_{n}}:=\{\left(  z,x,w\right)  \in\mathcal{F}%
^{n}|\;w_{i}=\infty\text{ with }1\leq i\leq n\ \ \text{implies}%
\]%
\[
x_{i}=-z-x_{1}-x_{2}-...-x_{i-1}\text{ and }x_{i+1}=...=x_{n}=0\}
\]
of $\mathcal{F}^{n}$, where $\sim$ is the equivalence relation generated by
\[
\left(  z,x,w\right)  \sim\left(  z,x,w_{1},...,w_{i}=\infty,\infty
,...,\infty\right)
\]
for all $\left(  z,x,w\right)  $ with $w_{i}=\infty$ for an $1\leq i\leq n$.
The unit space of $\mathfrak{F}_{n}$ is $Z:=\overline{{\mathbb{Z}}}_{\geq}%
^{n}/\sim$ where $\overline{{\mathbb{Z}}}_{\geq}^{n}$ is the unit space of
$\widetilde{\mathfrak{F}_{n}}\subset\mathcal{F}^{n}$ embedded in
$\widetilde{\mathfrak{F}_{n}}$ as the $\sim$-invariant subset $\left\{
0\right\}  \times\left\{  0\right\}  \times\overline{{\mathbb{Z}}}_{\geq}^{n}$.

Let $\pi_{n}$ denote the faithful *-representation of the groupoid C*-algebra
$C^{\ast}\left(  \mathfrak{F}_{n}\right)  $ canonically constructed on the
Hilbert space $\ell^{2}\left(  \mathbb{Z}\times\mathbb{Z}_{\geq}^{n}\right)
=\ell^{2}\left(  \mathbb{Z}\right)  \otimes\ell^{2}\left(  \mathbb{Z}_{\geq
}^{n}\right)  $ built from the open dense orbit $\mathbb{Z}_{\geq}^{n}$ in the
unit space $Z$ of $\mathfrak{F}_{n}$. For practical purposes, we often
identify $C^{\ast}\left(  \mathfrak{F}_{n}\right)  $ with the concrete
operator algebra $\pi_{n}\left(  C^{\ast}\left(  \mathfrak{F}_{n}\right)
\right)  $ without making explicit distinction. Note that by restricting
$\mathfrak{F}_{n}$ to the open subset $\mathbb{Z}_{\geq}^{n}$, we get the
groupoid $\mathfrak{F}_{n}|_{\mathbb{Z}_{\geq}^{n}}\cong\mathbb{Z}%
\times\left(  \left(  \mathbb{Z}^{n}\ltimes\mathbb{Z}^{n}\right)
|_{\mathbb{Z}_{\geq}^{n}}\right)  $ with
\[
C^{\ast}\left(  \mathfrak{F}_{n}|_{\mathbb{Z}_{\geq}^{n}}\right)  \cong
C^{\ast}\left(  \mathbb{Z}\right)  \otimes C^{\ast}\left(  \left(
\mathbb{Z}^{n}\ltimes\mathbb{Z}^{n}\right)  |_{\mathbb{Z}_{\geq}^{n}}\right)
\cong C\left(  \mathbb{T}\right)  \otimes\mathcal{K}\left(  \ell^{2}\left(
\mathbb{Z}_{\geq}^{n}\right)  \right)
\]
under the representation $\pi_{n}$, where $C^{\ast}\left(  \mathbb{Z}\right)
\cong C\left(  \mathbb{T}\right)  $ acts on $\ell^{2}\left(  \mathbb{Z}%
\right)  \cong L^{2}\left(  \mathbb{T}\right)  $ by multiplication operators,
and hence $C\left(  \mathbb{T}\right)  \otimes\mathcal{K}\left(  \ell
^{2}\left(  \mathbb{Z}_{\geq}^{n}\right)  \right)  $ can be viewed as a closed
ideal of $\pi_{n}\left(  C^{\ast}\left(  \mathfrak{F}_{n}\right)  \right)
\equiv C^{\ast}\left(  \mathfrak{F}_{n}\right)  $.

Note that the $\mathbb{Z}$-component of $\mathcal{F}^{n}\equiv{\mathbb{Z}%
}\times\left(  {\mathbb{Z}}^{n}\ltimes\overline{{\mathbb{Z}}}^{n}\right)
|_{\overline{{\mathbb{Z}}}_{\geq}^{n}}$ gives a grading on $C^{\ast}\left(
\mathfrak{F}_{n}\right)  $, decomposing it into (a completion of) a direct sum
of some subspaces index by $\mathbb{Z}$. More precisely,
$\widetilde{\mathfrak{F}_{n}}$ is the union of the pairwise disjoint closed
and open sets
\[
\left(  \widetilde{\mathfrak{F}_{n}}\right)  _{k}:=\left\{  \left(
k,x,w\right)  \in\widetilde{\mathfrak{F}_{n}}|\ \left(  x,w\right)
\in{\mathbb{Z}}^{n}\times\overline{{\mathbb{Z}}}_{\geq}^{n}\right\}
\]
with $k\in\mathbb{Z}$ which are invariant under the equivalence relation
$\sim$, so $\mathfrak{F}\equiv\widetilde{\mathfrak{F}_{n}}/\sim$ is the union
of the pairwise disjoint closed and open sets
\[
\left(  \mathfrak{F}_{n}\right)  _{k}:=\left(  \widetilde{\mathfrak{F}_{n}%
}\right)  _{k}/\sim
\]
and hence $C^{\ast}\left(  \mathfrak{F}_{n}\right)  $ is the closure of the
(algebraic) direct sum $\oplus_{k\in\mathbb{Z}}C_{c}\left(  \mathfrak{F}%
_{n}\right)  _{k}$ where $C_{c}\left(  \mathfrak{F}_{n}\right)  _{k}%
:=C_{c}\left(  \left(  \mathfrak{F}_{n}\right)  _{k}\right)  $. In fact, the
groupoid character $\left[  \left(  k,x,w\right)  \right]  \in\mathfrak{F}%
_{n}\mapsto t^{k}\in\mathbb{T}$ for any fixed $t\in\mathbb{T}\equiv U\left(
1\right)  $ defines an isometric *-automorphism of $L^{1}\left(
\mathfrak{F}_{n}\right)  $ and hence a C*-algebra automorphism $\rho\left(
t\right)  $ of $C^{\ast}\left(  \mathfrak{F}_{n}\right)  $, sending
$\delta_{\left[  \left(  k,x,w\right)  \right]  }$ to $t^{k}\delta_{\left[
\left(  k,x,w\right)  \right]  }$. Clearly $\rho:t\mapsto\rho\left(  t\right)
$ defines a $U\left(  1\right)  $-action on $C^{\ast}\left(  \mathfrak{F}%
_{n}\right)  $. The degree-$k$ spectral subspace $C^{\ast}\left(
\mathfrak{F}_{n}\right)  _{k}$ of $C^{\ast}\left(  \mathfrak{F}_{n}\right)  $
under the action $\rho$, i.e. the set consisting of all elements $a\in
C^{\ast}\left(  \mathfrak{F}_{n}\right)  $ with $\left(  \rho\left(  t\right)
\right)  \left(  a\right)  =t^{k}a$ for all $t\in\mathbb{T}$, is a closed
linear subspace of $C^{\ast}\left(  \mathfrak{F}_{n}\right)  $ containing
$C_{c}\left(  \mathfrak{F}_{n}\right)  _{k}$. Clearly $C^{\ast}\left(
\mathfrak{F}_{n}\right)  _{k}\cap C^{\ast}\left(  \mathfrak{F}_{n}\right)
_{k^{\prime}}=0$ for any $k\neq k^{\prime}$. The integration operator
\[
\Lambda_{k}:a\in C^{\ast}\left(  \mathfrak{F}_{n}\right)  \mapsto a_{k}%
:=\int_{\mathbb{T}}t^{-k}\left(  \rho\left(  t\right)  \right)  \left(
a\right)  dt\in C^{\ast}\left(  \mathfrak{F}_{n}\right)  _{k}\subset C^{\ast
}\left(  \mathfrak{F}_{n}\right)
\]
is a well-defined continuous projection onto $C^{\ast}\left(  \mathfrak{F}%
_{n}\right)  _{k}$ and eliminates $C^{\ast}\left(  \mathfrak{F}_{n}\right)
_{l}$ for all $l\neq k$, where $\mathbb{T}$ is endowed with the standard Haar
measure. Indeed for any $s\in\mathbb{T}$,%
\[
\left(  \rho\left(  s\right)  \right)  \left(  a_{k}\right)  =\int%
_{\mathbb{T}}t^{-k}\left(  \rho\left(  t\right)  \right)  \left(  \rho\left(
s\right)  a\right)  dt=s^{k}\int_{\mathbb{T}}\left(  st\right)  ^{-k}\left(
\rho\left(  st\right)  \right)  \left(  a\right)  dt
\]%
\[
=s^{k}\int_{\mathbb{T}}t^{-k}\left(  \rho\left(  t\right)  \right)  \left(
a\right)  dt=s^{k}a_{k},
\]
and for any $b\in C^{\ast}\left(  \mathfrak{F}_{n}\right)  _{l}$,
\[
\Lambda_{k}\left(  b\right)  =\int_{\mathbb{T}}t^{-k}\left(  \rho\left(
t\right)  \right)  \left(  b\right)  dt=\int_{\mathbb{T}}t^{-k}t^{l}%
bdt=\left(  \int_{\mathbb{T}}t^{l-k}dt\right)  b=\delta_{kl}b.
\]
So $\Lambda_{k}$'s are mutually orthogonal projections in the sense that
$\Lambda_{k}\circ\Lambda_{l}=\delta_{kl}\Lambda_{k}$. With the (algebraic) sum
$\sum_{k\in\mathbb{Z}}C_{c}\left(  \mathfrak{F}_{n}\right)  _{k}$ clearly
dense in $C^{\ast}\left(  \mathfrak{F}_{n}\right)  $ and $C_{c}\left(
\mathfrak{F}_{n}\right)  _{k}\subset C^{\ast}\left(  \mathfrak{F}_{n}\right)
_{k}$, we see that $\overline{C_{c}\left(  \mathfrak{F}_{n}\right)  _{k}%
}=C^{\ast}\left(  \mathfrak{F}_{n}\right)  _{k}$ by applying the projection
operator $\Lambda_{k}$ to any sequence in $\sum_{l\in\mathbb{Z}}C_{c}\left(
\mathfrak{F}_{n}\right)  _{l}$ converging to an element of $C^{\ast}\left(
\mathfrak{F}_{n}\right)  _{k}$. Furthermore we note that clearly $C_{c}\left(
\mathfrak{F}_{n}\right)  _{k}C_{c}\left(  \mathfrak{F}_{n}\right)  _{l}\subset
C_{c}\left(  \mathfrak{F}_{n}\right)  _{k+l}$ and $C^{\ast}\left(
\mathfrak{F}_{n}\right)  _{k}C^{\ast}\left(  \mathfrak{F}_{n}\right)
_{l}\subset C^{\ast}\left(  \mathfrak{F}_{n}\right)  _{k+l}$ for all
$k,l\in\mathbb{Z}$, i.e. $C_{c}\left(  \mathfrak{F}_{n}\right)  $ and
$C^{\ast}\left(  \mathfrak{F}_{n}\right)  $ are graded algebras (up to completion).

Recall that the group $U\left(  1\right)  \equiv\mathbb{T}$ acts on $C\left(
\mathbb{S}_{q}^{2n+1}\right)  $ by sending the standard generators
$u_{n+1,m}\in C\left(  SU_{q}\left(  n+1\right)  \right)  $, $1\leq m\leq
n+1$, of $C\left(  \mathbb{S}_{q}^{2n+1}\right)  $ to $tu_{n+1,m}$ for each
group element $t\in\mathbb{T}\subset\mathbb{C}$. This $U\left(  1\right)
$-action, denoted as $\tau_{t}$ for $t\in\mathbb{T}$, decomposes $C\left(
\mathbb{S}_{q}^{2n+1}\right)  $ into spectral subspaces $C\left(
\mathbb{S}_{q}^{2n+1}\right)  _{k}$ of degree $k\in\mathbb{Z}$ consisting of
elements $a\in C\left(  \mathbb{S}_{q}^{2n+1}\right)  $ satisfying $\tau
_{t}\left(  a\right)  =t^{k}a$ for all $t\in\mathbb{T}$. Each $u_{n+1,m}$ is
in the degree-$1$ spectral subspace $C\left(  \mathbb{S}_{q}^{2n+1}\right)
_{1}$. On the other hand, under the identification of $C\left(  \mathbb{S}%
_{q}^{2n+1}\right)  $ with $C^{\ast}\left(  \mathfrak{F}_{n}\right)  $
established in the work of \cite{Sh:cqg,Sh:qs}, each $u_{n+1,m}$ faithfully
represented as $t_{n+1}\otimes\gamma^{\otimes n+1-m}\otimes\alpha^{\ast
}\otimes1^{\otimes m-2}$ is identified with an element in $\overline
{C_{c}\left(  \mathfrak{F}_{n}\right)  _{1}}=C^{\ast}\left(  \mathfrak{F}%
_{n}\right)  _{1}$. So the grading on $C^{\ast}\left(  \mathfrak{F}%
_{n}\right)  $ by $C^{\ast}\left(  \mathfrak{F}_{n}\right)  _{k}$ coincides
with the grading on $C\left(  \mathbb{S}_{q}^{2n+1}\right)  $ by $C\left(
\mathbb{S}_{q}^{2n+1}\right)  _{k}$, i.e. $C\left(  \mathbb{S}_{q}%
^{2n+1}\right)  _{k}=C^{\ast}\left(  \mathfrak{F}_{n}\right)  _{k}$.

The degree-$0$ spectral subspace $C\left(  \mathbb{S}_{q}^{2n+1}\right)  _{0}%
$, or equivalently, the $U\left(  1\right)  $-invariant subalgebra $\left(
C\left(  \mathbb{S}_{q}^{2n+1}\right)  \right)  ^{U\left(  1\right)  }$ of
$C\left(  \mathbb{S}_{q}^{2n+1}\right)  $ can be naturally called the algebra
of quantum $\mathbb{C}P^{n}$, denoted as $C\left(  \mathbb{C}P_{q}^{n}\right)
$. The embedding $C\left(  \mathbb{C}P_{q}^{n}\right)  \subset C\left(
\mathbb{S}_{q}^{2n+1}\right)  \equiv C^{\ast}\left(  \mathfrak{F}_{n}\right)
$, or virtually the quantum quotient map $\mathbb{S}_{q}^{2n+1}\rightarrow
\mathbb{C}P_{q}^{n}$, is a quantum analogue of the Hopf principal $U\left(
1\right)  $-bundle $\mathbb{S}^{2n+1}\rightarrow\mathbb{C}P^{n}$. Furthermore
the degree-$k$ spectral subspaces $C\left(  \mathbb{S}_{q}^{2n+1}\right)
_{k}\equiv C^{\ast}\left(  \mathfrak{F}_{n}\right)  _{k}$ become the quantum
line bundles, denoted $L_{k}$, over $\mathbb{C}P_{q}^{n}$ associated with the
quantum principal $U\left(  1\right)  $-bundle $\mathbb{S}_{q}^{2n+1}%
\rightarrow\mathbb{C}P_{q}^{n}$. Note that in the context of groupoid
C*-algebras, $C\left(  \mathbb{C}P_{q}^{n}\right)  \equiv C\left(
\mathbb{S}_{q}^{2n+1}\right)  _{0}$ is canonically identified with the
groupoid C*-algebra $C^{\ast}\left(  \left(  \mathfrak{F}_{n}\right)
_{0}\right)  $ where $\left(  \mathfrak{F}_{n}\right)  _{0}$ is clearly an
open and closed subgroupoid of $\mathfrak{F}_{n}$. It is easy to see that the
unit space of $\left(  \mathfrak{F}_{n}\right)  _{0}\subset\mathfrak{F}_{n}$
is the same unit space $Z\equiv\overline{{\mathbb{Z}}}_{\geq}^{n}/\sim$ that
$\mathfrak{F}_{n}$ has.

On the other hand, the quantum complex projective space $U\left(  n\right)
_{q}\backslash SU\left(  n+1\right)  _{q}$ has been formulated and studied by
researchers from the viewpoint of quantum homogeneous space \cite{Na}. The
author showed in \cite{Sh:qcp} that such a quantum space can be concretely
realized by the C*-subalgebra generated by $u_{n+1,i}^{\ast}u_{n+1,j}$ with
$1\leq i,j\leq n+1$ in $C\left(  \mathbb{S}_{q}^{2n+1}\right)  $, and then
identified this C*-algebra with the groupoid C*-algebra $C^{\ast}%
(\mathfrak{T}_{n})$ of the subquotient groupoid $\mathfrak{T}_{n}%
:=\widetilde{\mathfrak{T}_{n}}/\sim$ of ${\mathbb{Z}}^{n}\ltimes
\overline{{\mathbb{Z}}}^{n}|_{\overline{{\mathbb{Z}}}_{\geq}^{n}}$, where%

\[
\widetilde{\mathfrak{T}_{n}}:=\{\left(  x,w\right)  \in{\mathbb{Z}}^{n}%
\ltimes\overline{{\mathbb{Z}}}^{n}|_{\overline{{\mathbb{Z}}}_{\geq}^{n}%
}:\;w_{i}=\infty\text{ with }1\leq i\leq n\ \ \text{implies}%
\]%
\[
x_{i}=-x_{1}-x_{2}-...-x_{i-1}\text{ and }x_{i+1}=...=x_{n}=0\}
\]
is a subgroupoid of ${\mathbb{Z}}^{n}\times\overline{{\mathbb{Z}}}%
^{n}|_{\overline{{\mathbb{Z}}}_{\geq}^{n}}$ and $\sim$ is the equivalence
relation generated by
\[
\left(  x,w\right)  \sim\left(  x,w_{1},...,w_{i}=\infty,\infty,...,\infty
\right)
\]
for all $(x,w)$ with $w_{i}=\infty$ for an $1\leq i\leq n$. It is easy to see
that $\left[  \left(  0,x,w\right)  \right]  \in\left(  \mathfrak{F}%
_{n}\right)  _{0}\mapsto\left[  \left(  x,w\right)  \right]  \in
\mathfrak{T}_{n}$ is a well-defined homeomorphic groupoid isomorphism, and
hence $C^{\ast}(\mathfrak{T}_{n})\cong C^{\ast}\left(  \left(  \mathfrak{F}%
_{n}\right)  _{0}\right)  $. So the quantum homogeneous space $U\left(
n\right)  _{q}\backslash SU\left(  n+1\right)  _{q}$ coincides with the
quantum complex projective space $\mathbb{C}P_{q}^{n}$ defined above, and the
results obtained in \cite{Sh:qcp} are valid for our study of the quantum
complex projective space $\mathbb{C}P_{q}^{n}$.

\section{Projections over $C\left(  \mathbb{S}_{q}^{2n+1}\right)  $}

In \cite{Sh:qcp}, taking the groupoid C*-algebra approach, we established an
inductive family of short exact sequences of C*-algebras
\[
0\rightarrow C\left(  \mathbb{T}\right)  \otimes\mathcal{K}\left(  \ell
^{2}\left(  \mathbb{Z}_{\geq}^{n}\right)  \right)  \rightarrow C\left(
\mathbb{S}_{q}^{2n+1}\right)  \rightarrow C\left(  \mathbb{S}_{q}%
^{2n-1}\right)  \rightarrow0.
\]
However for the purpose of classification of projections over $C\left(
\mathbb{S}_{q}^{2n+1}\right)  $, it turns out that another inductive family of
short exact sequences constructed below is more convenient.

Under the groupoid monomorphism
\[
\left(  z,x,w\right)  \in\widetilde{\mathfrak{F}_{n}}\mapsto\left(
z+x_{1},x,w\right)  \in\mathcal{F}^{n},
\]
$\widetilde{\mathfrak{F}_{n}}$ is mapped homeomorphically onto the image
$\widetilde{\mathfrak{F}_{n}}^{\prime}\subset\mathcal{F}^{n}$ consisting of
$\left(  z,x,w\right)  \in\mathcal{F}^{n}$ satisfying%
\[
\left\{
\begin{array}
[c]{lll}%
w_{1}=\infty\text{ }\ \Longrightarrow\ \text{\textquotedblleft}z=0\text{ and
}x_{2}=...=x_{n}=0\text{\textquotedblright}, &  & \\
w_{i}=\infty\text{ }\ \Longrightarrow\ \text{\textquotedblleft}x_{i}%
=-z-x_{2}-...-x_{i-1}\text{ and }x_{i+1}=...=x_{n}=0\text{\textquotedblright%
}, & \text{for } & 2\leq i\leq n,
\end{array}
\right.
\]
while the equivalence relation $\sim$ on $\widetilde{\mathfrak{F}_{n}}$
remains the same equivalence relation $\sim^{\prime}$ on
$\widetilde{\mathfrak{F}_{n}}^{\prime}$ that is generated by
\[
\left(  z,x,w\right)  \sim^{\prime}\left(  z,x,w_{1},...,w_{i}=\infty
,\infty,...,\infty\right)
\]
for all $\left(  z,x,w\right)  $ with $w_{i}=\infty$ for some $1\leq i\leq n$.
So we get a homeomorphic groupoid isomorphism%
\[
\gamma:\left[  \left(  z,x,w\right)  \right]  \in\widetilde{\mathfrak{F}_{n}%
}/\sim\equiv\mathfrak{F}_{n}\mapsto\left[  \left(  z+x_{1},x,w\right)
\right]  \in\widetilde{\mathfrak{F}_{n}}^{\prime}/\sim^{\prime}=:\mathfrak{F}%
_{n}^{\prime}.
\]
Note that the groupoid C*-algebra $C^{\ast}\left(  \mathfrak{F}_{n}^{\prime
}\right)  $ also has a faithful *-representation $\pi_{n}^{\prime}$
canonically constructed on the Hilbert space $\ell^{2}\left(  \mathbb{Z}%
\times\mathbb{Z}_{\geq}^{n}\right)  =\ell^{2}\left(  \mathbb{Z}\right)
\otimes\ell^{2}\left(  \mathbb{Z}_{\geq}^{n}\right)  $ built from the open
dense orbit $\mathbb{Z}_{\geq}^{n}$ in the unit space of $\mathfrak{F}%
_{n}^{\prime}$.

Let $m^{\left(  k\right)  }$ denote $\left(  m,...,m\right)  \in
\overline{{\mathbb{Z}}}^{k}$. Note that $\left(  \left\{  \infty\right\}
\times\overline{{\mathbb{Z}}}_{\geq}^{n-1}\right)  /\sim=\left\{  \left[
\infty^{\left(  n\right)  }\right]  \right\}  $ is a closed invariant subset
of the unit space $Z\equiv\overline{{\mathbb{Z}}}_{\geq}^{n}/\sim$ of
$\mathfrak{F}_{n}$ such that with a singleton unit space,
\[
\mathfrak{F}_{n}|_{\left\{  \left[  \infty^{\left(  n\right)  }\right]
\right\}  }\equiv\left\{  \left[  \left(  z,-z,0^{\left(  n-1\right)  }%
,\infty^{\left(  n\right)  }\right)  \right]  :z\in\mathbb{Z}\right\}
\cong\mathbb{Z}%
\]
as a group. On the other hand, the complement of $\left\{  \left[
\infty^{\left(  n\right)  }\right]  \right\}  $ in $Z$ is the open invariant
subset $O:=\left(  \mathbb{Z}_{\geq}\times\overline{{\mathbb{Z}}}_{\geq}%
^{n-1}\right)  /\sim$ such that $w_{1}\neq\infty$ for all $\left[  \left(
z,x,w\right)  \right]  \in\mathfrak{F}_{n}|_{O}$ and hence in $\gamma\left(
\left[  \left(  z,x,w\right)  \right]  \right)  =\left[  \left(
z+x_{1},x,w\right)  \right]  $, there is no non-trivial condition from the
definition of $\widetilde{\mathfrak{F}_{n}}^{\prime}/\sim^{\prime}$ imposed on
$\left(  x_{1},w_{1}\right)  $, while the non-trivial conditions from the
definition of $\widetilde{\mathfrak{F}_{n}}^{\prime}/\sim^{\prime}$ imposed on
the other components of $\gamma\left(  \left[  \left(  z,x,w\right)  \right]
\right)  $ $\ $match those in defining $\mathfrak{F}_{n-1}$. That is to say,
by rewriting $x_{1},w_{1}$ as the first two components of $\gamma\left(
\left[  \left(  z,x,w\right)  \right]  \right)  $, we have a homeomorphic
groupoid isomorphism from $\mathfrak{F}_{n}|_{O}$ onto the groupoid $\left(
\mathbb{Z}\ltimes\mathbb{Z}\right)  |_{\mathbb{Z}_{\geq}}\times\mathfrak{F}%
_{n-1}$, namely,%
\[
\gamma:\left[  \left(  z,x,w\right)  \right]  \in\mathfrak{F}_{n}|_{O}%
\mapsto\left(  x_{1},w_{1},\left[  \left(  z+x_{1},x_{2},..,x_{n}%
,w_{2},..,w_{n}\right)  \right]  \right)  \in\left(  \mathbb{Z}\ltimes
\mathbb{Z}\right)  |_{\mathbb{Z}_{\geq}}\times\mathfrak{F}_{n-1}%
\subset\mathfrak{F}_{n}^{\prime},
\]
which then induces a C*-algebra isomorphism
\[
\gamma_{\ast}:C^{\ast}\left(  \mathfrak{F}_{n}|_{O}\right)  \rightarrow
C^{\ast}\left(  \left(  \mathbb{Z}\ltimes\mathbb{Z}\right)  |_{\mathbb{Z}%
_{\geq}}\right)  \otimes C^{\ast}\left(  \mathfrak{F}_{n-1}\right)  \subset
C^{\ast}\left(  \mathfrak{F}_{n}^{\prime}\right)  .
\]

Note that $\pi_{n}^{\prime}=\pi_{0}\otimes\pi_{n-1}$ on $\gamma_{\ast}\left(
C^{\ast}\left(  \mathfrak{F}_{n}|_{O}\right)  \right)  $ where $\pi
_{0}:C^{\ast}\left(  \left(  \mathbb{Z}\ltimes\mathbb{Z}\right)
|_{\mathbb{Z}_{\geq}}\right)  \rightarrow\mathcal{K}\left(  \ell^{2}\left(
\mathbb{Z}_{\geq}\right)  \right)  $ is the well-known canonical faithful
representation, and the faithful representation
\[
\pi_{n}^{\prime}\circ\gamma_{\ast}:C^{\ast}\left(  \mathfrak{F}_{n}\right)
\rightarrow\mathcal{B}\left(  \ell^{2}\left(  \mathbb{Z}\times\mathbb{Z}%
_{\geq}^{n}\right)  \right)
\]
restricts to an isomorphism $C^{\ast}\left(  \mathfrak{F}_{n}|_{O}\right)
\cong\mathcal{K}\left(  \ell^{2}\left(  \mathbb{Z}_{\geq}\right)  \right)
\otimes C\left(  \mathbb{S}_{q}^{2n-1}\right)  $. So these invariant subsets
$\left\{  \left[  \infty^{n}\right]  \right\}  $ and $O$ give rise to a short
exact sequence%
\[
0\rightarrow\mathcal{K}\left(  \ell^{2}\left(  \mathbb{Z}_{\geq}\right)
\right)  \otimes C\left(  \mathbb{S}_{q}^{2n-1}\right)  \cong C^{\ast}\left(
\mathfrak{F}_{n}|_{O}\right)  \rightarrow C\left(  \mathbb{S}_{q}%
^{2n+1}\right)  \overset{\eta}{\rightarrow}C^{\ast}\left(  \mathfrak{F}%
_{n}|_{\left\{  \left[  \infty^{n}\right]  \right\}  }\right)  \cong C\left(
\mathbb{T}\right)  \rightarrow0.
\]

The set $T:=\left\{  \left[  \left(  z,x,w\right)  \right]  \in\mathfrak{F}%
_{n}:\ x_{1}=1=-z,\ x_{2}=...=x_{n}=0\right\}  $ is a compact open subset of
$\mathfrak{F}_{n}$, corresponding to the set $\left\{  \left[  \left(
0,x,w\right)  \right]  \in\mathfrak{F}_{n}^{\prime}:x=\left(  1,0^{\left(
n-1\right)  }\right)  \right\}  \subset\mathfrak{F}_{n}^{\prime}$ under the
isomorphism $\gamma$, and its characteristic function $\chi_{T}\in
C_{c}\left(  \mathfrak{F}_{n}\right)  \subset C^{\ast}\left(  \mathfrak{F}%
_{n}\right)  $ determines the operator
\[
\pi_{n}^{\prime}\left(  \gamma_{\ast}\left(  \chi_{T}\right)  \right)
=\mathcal{S}\otimes\mathrm{id}\in\mathcal{T}\otimes\pi_{n-1}\left(  C\left(
\mathbb{S}_{q}^{2n-1}\right)  \right)  \subset\mathcal{B}\left(  \ell
^{2}\left(  \mathbb{Z}_{\geq}\right)  \right)  \otimes\mathcal{B}\left(
\ell^{2}\left(  \mathbb{Z}\times\mathbb{Z}_{\geq}^{n-1}\right)  \right)
\]
where $\mathcal{S}$ is the unilateral shift operator generating the Toeplitz
algebra $\mathcal{T}$ with $\sigma\left(  \mathcal{S}\right)  =\mathrm{id}%
_{\mathbb{T}}$ for the symbol map $\sigma$ in the short exact sequence
$0\rightarrow\mathcal{K}\rightarrow\mathcal{T}\overset{\sigma}{\rightarrow
}C\left(  \mathbb{T}\right)  \rightarrow0$. Since the quotient map
$\eta:C^{\ast}\left(  \mathfrak{F}_{n}\right)  \rightarrow C^{\ast}\left(
\mathfrak{F}_{n}|_{\left\{  \left[  \infty^{n}\right]  \right\}  }\right)
\equiv C^{\ast}\left(  \mathbb{Z}\right)  $ restricts $\chi_{T}$ to
$\delta_{1}\in C_{c}\left(  \mathbb{Z}\right)  \equiv C_{c}\left(
\mathfrak{F}_{n}|_{\left\{  \left[  \infty^{n}\right]  \right\}  }\right)  $
yielding the function $\mathrm{id}_{\mathbb{T}}\in C\left(  \mathbb{T}\right)
\equiv C^{\ast}\left(  \mathbb{Z}\right)  $, we get
\[
C\left(  \mathbb{S}_{q}^{2n+1}\right)  \subset\mathcal{T}\otimes\pi
_{n-1}\left(  C\left(  \mathbb{S}_{q}^{2n-1}\right)  \right)  \equiv
\mathcal{T}\otimes C\left(  \mathbb{S}_{q}^{2n-1}\right)
\]
being the sum of $\mathcal{K}\otimes C\left(  \mathbb{S}_{q}^{2n-1}\right)  $
and $\mathcal{T}\otimes1_{C\left(  \mathbb{S}_{q}^{2n-1}\right)  }$, which
coincides with a description of $C\left(  \mathbb{S}_{q}^{2n+1}\right)  $ in
\cite{VaSo}.

The above surjective C*-algebra homomorphism $C\left(  \mathbb{S}_{q}%
^{2n+1}\right)  \overset{\eta}{\rightarrow}C\left(  \mathbb{T}\right)  $
facilitates the notion of rank for an equivalence class of idempotent $P\in
M_{\infty}\left(  C\left(  \mathbb{S}_{q}^{2n+1}\right)  \right)  $ over
$C\left(  \mathbb{S}_{q}^{2n+1}\right)  $, namely, the well-defined classical
rank of the vector bundle over $\mathbb{T}$ determined by the idempotent
$\eta\left(  P\right)  $ over $C\left(  \mathbb{T}\right)  $.

The set of equivalence classes of idempotents $P\in M_{\infty}\left(  C\left(
\mathbb{S}_{q}^{2n+1}\right)  \right)  $ equipped with the binary operation
$\boxplus$ becomes an abelian graded monoid
\[
\mathfrak{P}\left(  C\left(  \mathbb{S}_{q}^{2n+1}\right)  \right)
=\sqcup_{r=0}^{\infty}\mathfrak{P}_{r}\left(  C\left(  \mathbb{S}_{q}%
^{2n+1}\right)  \right)
\]
where $\mathfrak{P}_{r}\left(  C\left(  \mathbb{S}_{q}^{2n+1}\right)  \right)
$ is the set of all (equivalence classes of) idempotents over $C\left(
\mathbb{S}_{q}^{2n+1}\right)  $ of rank $r$, and
\[
\mathfrak{P}_{r}\left(  C\left(  \mathbb{S}_{q}^{2n+1}\right)  \right)
\boxplus\mathfrak{P}_{l}\left(  C\left(  \mathbb{S}_{q}^{2n+1}\right)
\right)  \subset\mathfrak{P}_{r+l}\left(  C\left(  \mathbb{S}_{q}%
^{2n+1}\right)  \right)
\]
for $r,l\geq0$. Clearly $\mathfrak{P}_{0}\left(  C\left(  \mathbb{S}%
_{q}^{2n+1}\right)  \right)  $ is a submonoid of $\mathfrak{P}\left(  C\left(
\mathbb{S}_{q}^{2n+1}\right)  \right)  $.

Now we can proceed to classify up to equivalence all projections over
$C\left(  \mathbb{S}_{q}^{2n+1}\right)  $ by induction on $n$, extending the
result obtained in \cite{Ba} for the case of $n=1$.

First we define some standard basic projections
\[
P_{j,k}:=\left\{
\begin{array}
[c]{lll}%
1_{\mathbb{T}}\otimes\left(  \left(  \otimes^{j-1}P_{1}\right)  \otimes
P_{k}\otimes\left(  \otimes^{n-j}I\right)  \right)  \in C\left(
\mathbb{T}\right)  \otimes\left(  \mathcal{K}^{+}\right)  ^{\otimes n}, &
\text{if } & k>0\text{ and }1\leq j\leq n\\
1_{\mathbb{T}}\otimes\left(  \boxplus^{k}I^{\otimes n}\right)  \equiv
1_{\mathbb{T}}\otimes\left(  \boxplus^{k}\left(  \otimes^{n}I\right)  \right)
\in M_{k}\left(  C\left(  \mathbb{T}\right)  \otimes\left(  \mathcal{K}%
^{+}\right)  ^{\otimes n}\right)  , & \text{if } & k\geq0\text{ and }j=0
\end{array}
\right.
\]
where $I$ stands for the unit of $\mathcal{K}^{+}$. Note that $P_{0,0}=0$.
(For the convenience of argument, we also use the symbol $P_{j,k}$ for the
case of $n=0$, by taking $\left(  \mathcal{K}^{+}\right)  ^{\otimes
0}:=\mathbb{C}$ and noting that $P_{0,k}=1_{\mathbb{T}}\otimes\left(
\boxplus^{k}1\right)  \in M_{k}\left(  \mathbb{C}\right)  $ for $k\geq0$ makes
sense, while $P_{j,k}$ with $1\leq j\leq n$ does not exist when $n=0$.)

We note that the basic projection $P_{j,k}$ with $j\geq1$ is implemented by
the characteristic function $\chi_{A_{j,k}}$ of the compact open subset
\[
A_{j,k}:=\left(  \left\{  0\right\}  \times\left\{  0\right\}  ^{\times
n}\times\left\{  0\right\}  ^{\times j-1}\times\left\{  0,1,..,k-1\right\}
\times\overline{{\mathbb{Z}}}_{\geq}^{n-j}\right)  /\sim
\]
of $\mathfrak{F}_{n}$ under both representations $\pi_{n}$ and $\pi
_{n}^{\prime}\circ\gamma_{\ast}$. So each $P_{j,k}$ with $j\geq1$ is a
projection in $C\left(  \mathbb{S}_{q}^{2n+1}\right)  $. On the other hand,
$P_{0,k}=\boxplus^{k}\tilde{I}$ is the identity projection in $M_{k}\left(
C\left(  \mathbb{S}_{q}^{2n+1}\right)  \right)  $, where $\tilde{I}$ is the
identity element of $C\left(  \mathbb{S}_{q}^{2n+1}\right)  $.

Recall that in the inductive family of short exact sequences%
\[
0\rightarrow C\left(  \mathbb{T}\right)  \otimes\mathcal{K}\left(  \ell
^{2}\left(  \mathbb{Z}_{\geq}^{n}\right)  \right)  \rightarrow C\left(
\mathbb{S}_{q}^{2n+1}\right)  \overset{\mu_{n}}{\rightarrow}C\left(
\mathbb{S}_{q}^{2n-1}\right)  \rightarrow0
\]
for $C\left(  \mathbb{S}_{q}^{2n+1}\right)  $ found in \cite{Sh:qcp}, the
quotient map $\mu_{n}:C\left(  \mathbb{S}_{q}^{2n+1}\right)  \rightarrow
C\left(  \mathbb{S}_{q}^{2n-1}\right)  $ is implemented by the restriction
map
\[
C^{\ast}\left(  \mathfrak{F}_{n}\right)  \rightarrow C^{\ast}\left(
\mathfrak{F}_{n}|_{\left(  \overline{{\mathbb{Z}}}_{\geq}^{n-1}\times\left\{
\infty\right\}  \right)  /\sim}\right)  \cong C^{\ast}\left(  \mathfrak{F}%
_{n-1}\right)  .
\]
For any $n\in\mathbb{N}$, a projection $P$ over $C\left(  \mathbb{S}%
_{q}^{2n+1}\right)  $ annihilated by $M_{\infty}\left(  \mu_{n}\right)  $ is a
projection in $M_{\infty}\left(  C\left(  \mathbb{T}\right)  \otimes
\mathcal{K}\left(  \ell^{2}\left(  \mathbb{Z}_{\geq}^{n}\right)  \right)
\right)  $ and hence has a well-defined finite operator-rank $d_{n}\left(
P\right)  \in\mathbb{Z}_{\geq}$, namely, the rank of the projection operator
$P\left(  t\right)  \in M_{\infty}\left(  \mathcal{K}\left(  \ell^{2}\left(
\mathbb{Z}_{\geq}^{n}\right)  \right)  \right)  $, independent of
$t\in\mathbb{T}$. If $P$ is not annihilated by $\mu_{n}$, then we assign
$d_{n}\left(  P\right)  :=\infty$. Note that $d_{n}\left(  P\right)  $ depends
only on the equivalence class of $P$ over $C\left(  \mathbb{S}_{q}%
^{2n+1}\right)  $. In the degenerate case of $n=0$, for a projection $P$ over
$C\left(  \mathbb{S}_{q}^{1}\right)  \equiv C\left(  \mathbb{T}\right)  $, we
define $d_{0}\left(  P\right)  $ to be the finite rank of projection $P\left(
t\right)  \in M_{\infty}\left(  \mathbb{C}\right)  $, independent of
$t\in\mathbb{T}$.

Now for a projection $P$ over $C\left(  \mathbb{S}_{q}^{2n+1}\right)  $, we
define for $0\leq l\leq n$,
\[
\rho_{l}\left(  P\right)  :=d_{l}\left(  \left(  \mu_{l+1}\circ\cdots\circ
\mu_{n-1}\circ\mu_{n}\right)  \left(  P\right)  \right)  \overset{\text{if
}l=n}{\equiv}d_{n}\left(  P\right)
\]
which depends only on the equivalence class of $P$ over $C\left(
\mathbb{S}_{q}^{2n+1}\right)  $ and gives us a well-defined monoid
homomorphism
\[
\rho_{l}:\left(  \mathfrak{P}\left(  C\left(  \mathbb{S}_{q}^{2n+1}\right)
\right)  ,\boxplus\right)  \rightarrow\left(  \mathbb{Z}_{\geq}\cup\left\{
\infty\right\}  ,+\right)  .
\]
It is easy to verify that
\[
\rho_{l}\left(  P_{j,k}\right)  =\left\{
\begin{array}
[c]{lll}%
\infty, & \text{if } & n-j>n-l\text{, i.e. }j<l\\
k, & \text{if } & j=l\\
0, & \text{if } & n-j<n-l\text{, i.e. }j>l
\end{array}
\right.
\]
which shows that these projections $P_{j,k}$ are mutually inequivalent over
$C\left(  \mathbb{S}_{q}^{2n+1}\right)  $ because $P_{j,k}$'s with different
indices $\left(  j,k\right)  $ are distinguished by the collection of
homomorphisms $\rho_{0},\rho_{1},...,\rho_{n}$.

\textbf{Theorem 1}. $\mathfrak{P}\left(  C\left(  \mathbb{S}_{q}%
^{2n+1}\right)  \right)  $ for $n\geq0$ is the disjoint union of
\[
\mathfrak{P}_{0}\left(  C\left(  \mathbb{S}_{q}^{2n+1}\right)  \right)
=\left\{  \left[  P_{0,0}\right]  \right\}  \cup\left\{  \left[
P_{j,k}\right]  :k>0\text{ and }1\leq j\leq n\right\}  ,
\]
containing pairwise distinct $\left[  P_{j,k}\right]  $'s indexed by $\left(
j,k\right)  $, and
\[
\mathfrak{P}_{k}\left(  C\left(  \mathbb{S}_{q}^{2n+1}\right)  \right)
=\left\{  \left[  P_{0,k}\right]  \right\}  \ \text{ a singleton\ for all
}k>0,
\]
and its monoid structure is explicitly determined by that
\[
\left[  P_{j,k}\right]  \boxplus\left[  P_{j^{\prime},k^{\prime}}\right]
=\left\{
\begin{array}
[c]{lll}%
\left[  P_{j,k}\right]  , & \text{if } & 0\leq j<j^{\prime}\text{ and
}k,k^{\prime}>0,\\
\left[  P_{j,k+k^{\prime}}\right]  , & \text{if } & j=j^{\prime}\geq0.
\end{array}
\right.
\]
So $\left[  P_{j,k}\right]  =0$ in $K_{0}\left(  C\left(  \mathbb{S}%
_{q}^{2n+1}\right)  \right)  $ if and only if $1\leq j\leq n$ or $j=k=0$.

\textbf{Proof}. Knowing that $P_{j,k}$ are mutually inequivalent, we only need
to show that any projection over $C\left(  \mathbb{S}_{q}^{2n+1}\right)  $ is
equivalent to one of these $P_{j,k}$'s and verify the stated monoid structure.

We prove by induction on $n\geq0$.

When $n=0$, $C\left(  \mathbb{S}_{q}^{2n+1}\right)  =C\left(  \mathbb{T}%
\right)  $ and it is well known from algebraic topology about vector bundles
over $\mathbb{T}$ that isomorphism classes of (complex) vector bundles over
$\mathbb{T}$ are faithfully represented by trivial vector bundles, i.e.
$\mathfrak{P}_{0}\left(  C\left(  \mathbb{T}\right)  \right)  =\left\{
0\right\}  \equiv\left\{  \left[  P_{0,0}\right]  \right\}  $ while
$\mathfrak{P}_{k}\left(  C\left(  \mathbb{T}\right)  \right)  =\left\{
\left[  P_{0,k}\right]  \right\}  $\ for $k>0$. Then the statements of this
theorem for $n=0$ are clearly verified.

Now assume that the statements hold for $C\left(  \mathbb{S}_{q}%
^{2n-1}\right)  $. We need to show that they also hold for $C\left(
\mathbb{S}_{q}^{2n+1}\right)  $.

Since any complex vector bundle over $\mathbb{T}$ is trivial, any idempotent
over $C\left(  \mathbb{T}\right)  $ is equivalent to the standard projection
$1\otimes P_{m}\in C\left(  \mathbb{T}\right)  \otimes M_{\infty}\left(
\mathbb{C}\right)  $ for some $m\in\mathbb{Z}_{\geq}$. So for any nonzero
idempotent $P\in M_{\infty}\left(  C\left(  \mathbb{S}_{q}^{2n+1}\right)
\right)  $ over $C\left(  \mathbb{S}_{q}^{2n+1}\right)  $, there is some $U\in
GL_{\infty}\left(  C\left(  \mathbb{T}\right)  \right)  $ such that
\[
U\eta\left(  P\right)  U^{-1}=1\otimes P_{m}=\eta\left(  \boxplus^{m}\tilde
{I}\right)
\]
for some $m\in\mathbb{Z}_{\geq}$ where $\tilde{I}$ is the identity element of
$C\left(  \mathbb{S}_{q}^{2n+1}\right)  $ viewed as the identity element in
$\left(  \mathcal{K}\otimes C\left(  \mathbb{S}_{q}^{2n-1}\right)  \right)
^{+}\subset C\left(  \mathbb{S}_{q}^{2n+1}\right)  $, and hence $VPV^{-1}%
-\boxplus^{m}\tilde{I}\in M_{\infty}\left(  \mathcal{K}\otimes C\left(
\mathbb{S}_{q}^{2n-1}\right)  \right)  $ for any lift $V\in GL_{\infty}\left(
C\left(  \mathbb{S}_{q}^{2n+1}\right)  \right)  $ (which exists) of $U\boxplus
U^{-1}\in GL_{\infty}^{0}\left(  C\left(  \mathbb{T}\right)  \right)  $ along
$\eta$. Replacing $P$ by the equivalent $VPV^{-1}$, we may assume that
$P\in\left(  \boxplus^{m}\tilde{I}\right)  +M_{r-1}\left(  \mathcal{K}\otimes
C\left(  \mathbb{S}_{q}^{2n-1}\right)  \right)  $ for some large $r\geq m+1$.
Now since $M_{\infty}\left(  \mathbb{C}\right)  $ is dense in $\mathcal{K}$,
there is an idempotent $Q\in\left(  \boxplus^{m}\tilde{I}\right)
+M_{r-1}\left(  M_{N-1}\left(  \mathbb{C}\right)  \otimes C\left(
\mathbb{S}_{q}^{2n-1}\right)  \right)  $ sufficiently close to and hence
equivalent to $P$ over $C\left(  \mathbb{S}_{q}^{2n+1}\right)  $ for some
large $N$. So replacing $P$ by $Q$, we may assume that
\[
K:=P-\boxplus^{m}\tilde{I}\in M_{r-1}\left(  M_{N-1}\left(  \mathbb{C}\right)
\otimes C\left(  \mathbb{S}_{q}^{2n-1}\right)  \right)  .
\]

Rearranging the entries of $P\equiv K+\boxplus^{m}\tilde{I}\in M_{r-1}\left(
C\left(  \mathbb{S}_{q}^{2n+1}\right)  \right)  \subset M_{r}\left(  C\left(
\mathbb{S}_{q}^{2n+1}\right)  \right)  $ via conjugation by the unitary
\[
U_{r,N}:=\sum_{j=1}^{r-1}\left[  e_{jj}\otimes\left(  \left(  \mathcal{S}%
\otimes\mathrm{id}\right)  ^{\ast}\right)  ^{N}+e_{rj}\otimes\left(  \left(
\mathcal{S}\otimes\mathrm{id}\right)  ^{\left(  j-1\right)  N}P_{N}\right)
\right]  +e_{rr}\otimes\left(  \mathcal{S}\otimes\mathrm{id}\right)  ^{\left(
r-1\right)  N}%
\]%
\[
\in M_{r}\left(  \mathbb{C}\right)  \otimes C\left(  \mathbb{S}_{q}%
^{2n+1}\right)  \equiv M_{r}\left(  C\left(  \mathbb{S}_{q}^{2n+1}\right)
\right)
\]
we get the idempotent%
\[
U_{r,N}PU_{r,N}^{-1}\equiv U_{r,N}\left(  P\boxplus0\right)  U_{r,N}%
^{-1}=\left(  \left(  \boxplus^{m}\tilde{I}\right)  \boxplus\left(
\boxplus^{r-1-m}0\right)  \right)  \boxplus R
\]
for some idempotent
\[
R\in M_{\left(  r-1\right)  N}\left(  C\left(  \mathbb{S}_{q}^{2n-1}\right)
\right)  \subset\mathcal{K}\otimes C\left(  \mathbb{S}_{q}^{2n-1}\right)
\subset C\left(  \mathbb{S}_{q}^{2n+1}\right)
\]
which has rank at least $m$ as an idempotent over $C\left(  \mathbb{S}%
_{q}^{2n-1}\right)  $ since it contains $m$ copies of the identity element
$\tilde{I}^{\prime}$ of $C\left(  \mathbb{S}_{q}^{2n-1}\right)  $ as diagonal
$\boxplus$-summands, relocated from the $N$-th diagonal entry of each of the
$m$ copies of $\tilde{I}$ in $P$.

Now by the induction hypothesis, the idempotent $R\in M_{\left(  r-1\right)
N}\left(  C\left(  \mathbb{S}_{q}^{2n-1}\right)  \right)  $ is equivalent over
$C\left(  \mathbb{S}_{q}^{2n-1}\right)  $ to $P_{j,k}^{\prime}$ (denoting a
standard projection $P_{j,k}$ over $C\left(  \mathbb{S}_{q}^{2n-1}\right)  $)
with $0\leq j\leq n-1$ and $k>0$, which is identified with%
\[
P_{j+1,k}\equiv\left\{
\begin{array}
[c]{lll}%
P_{1}\otimes P_{j,k}^{\prime}\in P_{1}\otimes C\left(  \mathbb{S}_{q}%
^{2n-1}\right)  \subset\mathcal{K}\otimes C\left(  \mathbb{S}_{q}%
^{2n-1}\right)  \subset C\left(  \mathbb{S}_{q}^{2n+1}\right)  , & \text{if }
& j>0\\
P_{k}\otimes\tilde{I}^{\prime}\in P_{k}\otimes C\left(  \mathbb{S}_{q}%
^{2n-1}\right)  \subset\mathcal{K}\otimes C\left(  \mathbb{S}_{q}%
^{2n-1}\right)  \subset C\left(  \mathbb{S}_{q}^{2n+1}\right)  , & \text{if }
& j=0
\end{array}
\right.  ,
\]
i.e. $WRW^{-1}=P_{j,k}^{\prime}\equiv P_{j+1,k}$ for some invertible $W\in
M_{N^{\prime}}\left(  C\left(  \mathbb{S}_{q}^{2n-1}\right)  \right)  $ with
$N^{\prime}\geq\left(  r-1\right)  N$. Note that if $m>0$, then $R$ has a
positive rank as an idempotent over $C\left(  \mathbb{S}_{q}^{2n-1}\right)  $
and hence $j=0$. Since
\[
W+\left(  \tilde{I}-P_{N^{\prime}}\otimes\tilde{I}^{\prime}\right)  \in\left(
\mathcal{K}\otimes C\left(  \mathbb{S}_{q}^{2n-1}\right)  \right)  ^{+}\subset
C\left(  \mathbb{S}_{q}^{2n+1}\right)  ,
\]
we get $\left(  \left(  \boxplus^{m}\tilde{I}\right)  \boxplus\left(
\boxplus^{r-1-m}0\right)  \right)  \boxplus R$ equivalent over $C\left(
\mathbb{S}_{q}^{2n+1}\right)  $ to the projection
\[
\left(  \left(  \boxplus^{m}\tilde{I}\right)  \boxplus\left(  \boxplus
^{r-1-m}0\right)  \right)  \boxplus P_{j+1,k}%
\]
where $j=0$ if $m>0$.

If $m=0$, then clearly $P$ is equivalent over $C\left(  \mathbb{S}_{q}%
^{2n+1}\right)  $ to $P_{j+1,k}\in\mathcal{K}\otimes C\left(  \mathbb{S}%
_{q}^{2n-1}\right)  $ with $j+1>0$ and hence is of rank $0$. (We assumed
$P\ $nonzero, so $k>0$.)

If $m\in\mathbb{N}$ and hence $j=0$, then $P_{j+1,k}=P_{1,k}\equiv
P_{k}\otimes\tilde{I}^{\prime}$ and we can rearrange entries of $\left(
\left(  \boxplus^{m}\tilde{I}\right)  \boxplus\left(  \boxplus^{r-1-m}%
0\right)  \right)  \boxplus P_{1,k}$ via conjugation by the unitary%
\[
U_{l}:=e_{11}\otimes\left(  \mathcal{S}^{k}\otimes\mathrm{id}\right)
+e_{1r}\otimes P_{k}+\sum_{j=2}^{r-1}e_{jj}\otimes\tilde{I}+e_{rr}%
\otimes\left(  \mathcal{S}^{k}\otimes\mathrm{id}\right)  ^{\ast}%
\]%
\[
\in M_{r}\left(  \mathbb{C}\right)  \otimes C\left(  \mathbb{S}_{q}%
^{2n+1}\right)  \equiv M_{r}\left(  C\left(  \mathbb{S}_{q}^{2n+1}\right)
\right)
\]
to get $P$ equivalent over $C\left(  \mathbb{S}_{q}^{2n+1}\right)  $ to
\[
U_{l}\left(  \left(  \left(  \boxplus^{m}\tilde{I}\right)  \boxplus\left(
\boxplus^{r-1-m}0\right)  \right)  \boxplus P_{1,k}\right)  U_{l}^{-1}=\left(
\boxplus^{m}\tilde{I}\right)  \boxplus\left(  \boxplus^{r-m}0\right)
\equiv\boxplus^{m}\tilde{I}\equiv P_{0,m}%
\]
which is of rank $m\in\mathbb{N}$.

So we have proved the description of the sets $\mathfrak{P}_{k}\left(
C\left(  \mathbb{S}_{q}^{2n+1}\right)  \right)  $ in the theorem. It remains
to verify the monoid structure of $\mathfrak{P}\left(  C\left(  \mathbb{S}%
_{q}^{2n+1}\right)  \right)  $.

By specializing the above analysis for $P\equiv K+\boxplus^{m}\tilde{I}$ to
the case of $K=\left(  \boxplus^{m}0\right)  \boxplus P_{j+1,k}$, we have
already established that
\[
P_{0,m}\boxplus P_{j+1,k}\equiv\left(  \boxplus^{m}\tilde{I}\right)  \boxplus
P_{j+1,k}\sim_{C\left(  \mathbb{S}_{q}^{2n+1}\right)  }P_{0,m}%
\]
for all $m\in\mathbb{N}$ and $j+1>0$, while $\left[  P_{0,k}\right]
\boxplus\left[  P_{0,k^{\prime}}\right]  =\left[  P_{0,k+k^{\prime}}\right]  $
is obvious.

On the other hand, by induction hypothesis,
\[
P_{j,k}^{\prime}\boxplus P_{j^{\prime},k^{\prime}}^{\prime}\sim_{C\left(
\mathbb{S}_{q}^{2n-1}\right)  }\left\{
\begin{array}
[c]{lll}%
P_{j,k}^{\prime}, & \text{if } & 0\leq j<j^{\prime}\\
P_{j,k+k^{\prime}}^{\prime}, & \text{if } & j=j^{\prime}\geq0
\end{array}
\right.  .
\]
Now by applying $P_{1}\otimes\cdot$ to both sides of this equivalence, we get%
\[
P_{j+1,k}\boxplus P_{j^{\prime}+1,k^{\prime}}\sim_{C\left(  \mathbb{S}%
_{q}^{2n+1}\right)  }\left\{
\begin{array}
[c]{lll}%
P_{j+1,k}, & \text{if } & 1\leq j+1<j^{\prime}+1\\
P_{j+1,k+k^{\prime}}, & \text{if } & j+1=j^{\prime}+1\geq1
\end{array}
\right.
\]
since if an invertible $U\in M_{N}\left(  C\left(  \mathbb{S}_{q}%
^{2n-1}\right)  \right)  $ with $N$ sufficiently large conjugates an
idempotent $P$ over $C\left(  \mathbb{S}_{q}^{2n-1}\right)  $ to an idempotent
$Q$, then
\[
\left(  P_{1}\otimes U_{ij}\right)  _{i,j=1}^{N}+\boxplus^{N}\left(  \tilde
{I}-P_{1}\otimes\tilde{I}^{\prime}\right)  \in M_{N}\left(  C\left(
\mathbb{S}_{q}^{2n+1}\right)  \right)
\]
is an invertible conjugating $P_{1}\otimes P$ to $P_{1}\otimes Q$.

Now we have established all the monoid structure rules for $\mathfrak{P}%
\left(  C\left(  \mathbb{S}_{q}^{2n+1}\right)  \right)  $.

$\square$

\textbf{Remark}. The last part of the above proof about the monoid structure
of $\mathfrak{P}\left(  C\left(  \mathbb{S}_{q}^{2n+1}\right)  \right)  $ can
be avoided by applying the injective monoid homomorphism $\rho$ of the
following Corollary 2 to both sides of the equivalence relations describing
the monoid structure.

\textbf{Corollary 1}. All projections over $C\left(  \mathbb{S}_{q}%
^{2n+1}\right)  $ of strictly positive rank are trivial. The cancellation law
holds for projections of rank $\geq1$, but fails for projections of rank $0$
in case of $n>0$.

\textbf{Proof}. The only equivalence class of projection of a fixed rank $k>0$
is the trivial projection $\left[  P_{0,k}\right]  =\left[  \boxplus^{k}%
\tilde{I}\right]  $ classified above. By counting the rank, it is clear that
if $\boxplus^{k}\tilde{I}$ and $\boxplus^{k^{\prime}}\tilde{I}$ are stably
equivalent, then $k=k^{\prime}$. So the cancellation law holds for projections
of rank $\geq1$.

On the other hand, since for any distinct pairs $\left(  j,k\right)  $ and
$\left(  j^{\prime},k^{\prime}\right)  $ with $1\leq j,j^{\prime}\leq n$ and
$k,k^{\prime}>0$, $\left[  P_{j,k}\right]  \neq\left[  P_{j^{\prime}%
,k^{\prime}}\right]  $ but
\[
\left[  P_{j,k}\right]  \boxplus\left[  P_{0,1}\right]  =\left[
P_{0,1}\right]  =\left[  P_{j^{\prime},k^{\prime}}\right]  \boxplus\left[
P_{0,1}\right]  ,
\]
the cancellation law fails for such rank-$0$ projections $P_{j,k}$ and
$P_{j^{\prime},k^{\prime}}$.

$\square$

\textbf{Corollary 2}. The monoid $\mathfrak{P}\left(  C\left(  \mathbb{S}%
_{q}^{2n+1}\right)  \right)  $ is a submonoid of $\prod_{0\leq l\leq
n}\overline{\mathbb{Z}_{\geq}}$ via the injective monoid homomorphism
\[
\rho:P\in\mathfrak{P}\left(  C\left(  \mathbb{S}_{q}^{2n+1}\right)  \right)
\mapsto\prod_{0\leq l\leq n}\rho_{l}\left(  P\right)  \in\prod_{0\leq l\leq
n}\overline{\mathbb{Z}_{\geq}}.
\]

\textbf{Proof}. $\rho$ is injective since we already know that $\rho_{l}$'s
can distinguish the standard projections $P_{j,k}$ which have been shown to
constitute the whole monoid $\mathfrak{P}\left(  C\left(  \mathbb{S}%
_{q}^{2n+1}\right)  \right)  $.

$\square$

\section{Generating Projections of $K_{0}\left(  C\left(  \mathbb{C}P_{q}%
^{n}\right)  \right)  $}

In this section, we present a set of elementary projections over $C\left(
\mathbb{C}P_{q}^{n}\right)  $, whose $K_{0}$-classes form a set of free
generators of the abelian group $K_{0}\left(  C\left(  \mathbb{C}P_{q}%
^{n}\right)  \right)  $. We remark that a fascinating geometric construction
of free generators of $K_{0}\left(  C\left(  \mathbb{C}P_{q}^{n}\right)
\right)  $ has been found by D'Andrea and Landi in \cite{dAnLa}.

As discussed before, by restricting to the degree-$0$ part of the groupoid
$\mathfrak{F}_{n}$ consisting of exactly those $\left[  \left(  z,x,w\right)
\right]  $ with $z=0$, we get a subgroupoid $\left(  \mathfrak{F}_{n}\right)
_{0}$ which realizes $C\left(  \mathbb{C}P_{q}^{n}\right)  $ as a groupoid
C*-algebra $C^{\ast}\left(  \left(  \mathfrak{F}_{n}\right)  _{0}\right)  $.
Roughly speaking, $\left(  \mathfrak{F}_{n}\right)  _{0}$ can be extracted
from $\mathfrak{F}_{n}$ by simply ignoring or removing the $z$-component of
the elements $\left[  \left(  z,x,w\right)  \right]  $. Note that if $\left[
\left(  0,x,w\right)  \right]  \in\left(  \mathfrak{F}_{n}\right)  _{0}$ and
$w_{1}=\infty$, then $x_{1}=0$ by the defining condition on $\mathfrak{F}_{n}%
$. Furthermore since clearly $\pi_{n}\left(  C^{\ast}\left(  \left(
\mathfrak{F}_{n}\right)  _{0}\right)  \right)  \subset\mathrm{id}_{\ell
^{2}\left(  \mathbb{Z}\right)  }\otimes\mathcal{B}\left(  \ell^{2}\left(
\mathbb{Z}_{\geq}^{n}\right)  \right)  $, we will ignore the factor
$\mathrm{id}_{\ell^{2}\left(  \mathbb{Z}\right)  }$ and view $\pi
_{n}|_{C^{\ast}\left(  \left(  \mathfrak{F}_{n}\right)  _{0}\right)  }$ as a
faithful representation of $C^{\ast}\left(  \left(  \mathfrak{F}_{n}\right)
_{0}\right)  $ on $\ell^{2}\left(  \mathbb{Z}_{\geq}^{n}\right)  $ instead of
on $\ell^{2}\left(  \mathbb{Z}\times\mathbb{Z}_{\geq}^{n}\right)  $.

In \cite{Sh:qcp}, by considering the closed invariant subset $\left(
\overline{\mathbb{Z}}_{\geq}^{n-1}\times\left\{  \infty\right\}  \right)
/\sim$ (i.e. $\left\{  \left[  w\right]  :w\in\overline{\mathbb{Z}}_{\geq
}^{n-1}\times\left\{  \infty\right\}  \right\}  $ even though $\overline
{\mathbb{Z}}_{\geq}^{n-1}\times\left\{  \infty\right\}  $ is not really $\sim
$-invariant in the unit space $\overline{{\mathbb{Z}}}_{\geq}^{n}$ of
$\widetilde{\mathfrak{F}_{n}}\subset\mathcal{F}^{n}$) and its complement
$O_{0}$ in the unit space $Z$ of $\left(  \mathfrak{F}_{n}\right)  _{0}$ (and
of $\mathfrak{F}_{n}$ as well), we get the following short exact sequence%
\[
0\rightarrow\mathcal{K}\left(  \ell^{2}\left(  \mathbb{Z}_{\geq}^{n}\right)
\right)  \cong C^{\ast}\left(  \left(  \mathfrak{F}_{n}\right)  _{0}|_{O_{0}%
}\right)  \rightarrow C\left(  \mathbb{C}P_{q}^{n}\right)  \overset{\nu
}{\rightarrow}C^{\ast}\left(  \mathfrak{F}_{n}|_{\left(  \overline{\mathbb{Z}%
}_{\geq}^{n-1}\times\left\{  \infty\right\}  \right)  /\sim}\right)  \cong
C\left(  \mathbb{C}P_{q}^{n-1}\right)  \rightarrow0
\]
with $\left(  \mathfrak{F}_{n}\right)  _{0}|_{O_{0}}\cong\left(
\mathbb{Z}^{n}\ltimes\mathbb{Z}^{n}\right)  |_{\mathbb{Z}_{\geq}^{n}}$. Thus
we get the following 6-term exact sequence
\[%
\begin{array}
[c]{ccccccc}%
\mathbb{Z}= & K_{0}\left(  \mathcal{K}\left(  \ell^{2}\left(  \mathbb{Z}%
_{\geq}^{n}\right)  \right)  \right)  & \rightarrow & K_{0}\left(  C\left(
\mathbb{C}P_{q}^{n}\right)  \right)  & \overset{\nu_{\ast}}{\rightarrow} &
K_{0}\left(  C\left(  \mathbb{C}P_{q}^{n-1}\right)  \right)  & \\
& \uparrow &  &  &  & \downarrow & \\
& K_{1}\left(  C\left(  \mathbb{C}P_{q}^{n-1}\right)  \right)  & \leftarrow &
K_{1}\left(  C\left(  \mathbb{C}P_{q}^{n}\right)  \right)  & \leftarrow &
K_{1}\left(  \mathcal{K}\left(  \ell^{2}\left(  \mathbb{Z}_{\geq}^{n}\right)
\right)  \right)  & =0.
\end{array}
\]

By an induction on $n\geq1$, we can establish $K_{0}\left(  C\left(
\mathbb{C}P_{q}^{n}\right)  \right)  \cong\mathbb{Z}^{n+1}$ and $K_{1}\left(
C\left(  \mathbb{C}P_{q}^{n}\right)  \right)  =0$. In fact, in the case of
$n=1$, we have $K_{i}\left(  C\left(  \mathbb{C}P_{q}^{0}\right)  \right)
=K_{i}\left(  \mathbb{C}\right)  \cong\delta_{0i}\mathbb{Z}$ and hence
$K_{0}\left(  C\left(  \mathbb{C}P_{q}^{1}\right)  \right)  \cong%
\mathbb{Z}\oplus\mathbb{Z}$ and $K_{1}\left(  C\left(  \mathbb{C}P_{q}%
^{1}\right)  \right)  =0$. For $n>1$, the induction hypothesis $K_{1}\left(
C\left(  \mathbb{C}P_{q}^{n-1}\right)  \right)  =0$ and $K_{0}\left(  C\left(
\mathbb{C}P_{q}^{n-1}\right)  \right)  \cong\mathbb{Z}^{n}$ forces
\[
K_{0}\left(  C\left(  \mathbb{C}P_{q}^{n}\right)  \right)  \cong K_{0}\left(
\mathcal{K}\right)  \oplus K_{0}\left(  C\left(  \mathbb{C}P_{q}^{n-1}\right)
\right)  \equiv\mathbb{Z}\oplus\mathbb{Z}^{n}=\mathbb{Z}^{n+1}%
\]
and also $K_{1}\left(  C\left(  \mathbb{C}P_{q}^{n}\right)  \right)  =0$ in
the above 6-term exact sequence.

The above induction can be refined to get the following stronger result. First
we note that
\[
P_{j,k}\equiv1_{\mathbb{T}}\otimes\left(  \left(  \otimes^{j-1}P_{1}\right)
\otimes P_{k}\otimes\left(  \otimes^{n-j}I\right)  \right)  \in C\left(
\mathbb{S}_{q}^{2n+1}\right)  \subset C\left(  \mathbb{T}\right)
\otimes\mathcal{B}\left(  \ell^{2}\left(  \mathbb{Z}_{\geq}^{n}\right)
\right)
\]
with $0<j\leq n$ is a projection in $C\left(  \mathbb{C}P_{q}^{n}\right)
\subset C\left(  \mathbb{S}_{q}^{2n+1}\right)  $, and can be identified with
\[
\left(  \left(  \otimes^{j-1}P_{1}\right)  \otimes P_{k}\otimes\left(
\otimes^{n-j}I\right)  \right)  \in C\left(  \mathbb{C}P_{q}^{n}\right)
\subset\mathcal{B}\left(  \ell^{2}\left(  \mathbb{Z}_{\geq}^{n}\right)
\right)  .
\]
From now on, we view $P_{j,k}$ with $0<j\leq n$ as the latter elementary
tensor product lying in $C\left(  \mathbb{C}P_{q}^{n}\right)  $. On the other
hand, clearly the trivial projection $P_{0,k}$ of rank $k$ over $C\left(
\mathbb{S}_{q}^{2n+1}\right)  $ is also a trivial projection of rank $k$ over
$C\left(  \mathbb{C}P_{q}^{n}\right)  $.

\textbf{Theorem 2}. The standard projections $P_{j,1}\equiv\left(  \otimes
^{j}P_{1}\right)  \otimes\left(  \otimes^{n-j}I\right)  $ over $C\left(
\mathbb{C}P_{q}^{n}\right)  $ with $0\leq j\leq n$ are inequivalent over
$C\left(  \mathbb{C}P_{q}^{n}\right)  $ and their equivalence classes (over
$C\left(  \mathbb{C}P_{q}^{n}\right)  $, not over $C\left(  \mathbb{S}%
_{q}^{2n+1}\right)  $) form a set of free generators of $K_{0}\left(  C\left(
\mathbb{C}P_{q}^{n}\right)  \right)  \cong\mathbb{Z}^{n+1}$.

\textbf{Proof}. Since $P_{j,1}$ are inequivalent over $C\left(  \mathbb{S}%
_{q}^{2n+1}\right)  $, they are clearly inequivalent over the subalgebra
$C\left(  \mathbb{C}P_{q}^{n}\right)  $. Now we prove by induction on $n\geq1$
that $\left[  P_{j,1}\right]  $ with $0\leq j\leq n$ form a set of free
generators of $K_{0}\left(  C\left(  \mathbb{C}P_{q}^{n}\right)  \right)  $.

For $n=1$, it is well-known that $\mathcal{K}\left(  \ell^{2}\left(
\mathbb{Z}_{\geq}\right)  \right)  ^{+}\cong C\left(  \mathbb{C}P_{q}%
^{1}\right)  $ has $\left[  P_{1}\right]  \equiv\left[  P_{1,1}\right]  $ and
$\left[  I\right]  \equiv\left[  P_{0,1}\right]  $ as free generators of its
$K_{0}$-group $K_{0}\left(  \mathcal{K}\left(  \ell^{2}\left(  \mathbb{Z}%
_{\geq}\right)  \right)  ^{+}\right)  \cong\mathbb{Z}^{2}$.

For $n>1$, $K_{0}\left(  \mathcal{K}\left(  \ell^{2}\left(  \mathbb{Z}_{\geq
}^{n}\right)  \right)  \right)  \cong\mathbb{Z}$ has $\left[  \otimes^{n}%
P_{1}\right]  \equiv\left[  P_{n,1}\right]  $ as a free generator, while by
induction hypothesis, $K_{0}\left(  C\left(  \mathbb{C}P_{q}^{n-1}\right)
\right)  \cong\mathbb{Z}^{n}$ has $\left[  P_{j,1}^{\prime}\right]
\equiv\left[  \left(  \otimes^{j}P_{1}\right)  \otimes\left(  \otimes
^{n-1-j}I\right)  \right]  $ with $0\leq j\leq n-1$ as free generators. Now
with $\nu_{\ast}\left(  \left[  P_{j,1}\right]  \right)  \equiv\nu_{\ast
}\left(  \left[  P_{j,1}^{\prime}\otimes I\right]  \right)  =\left[
P_{j,1}^{\prime}\right]  $ for all $0\leq j\leq n-1$, it is easy to see from
the above 6-term exact sequence that $\left[  P_{j,1}\right]  $ for $0\leq
j\leq n-1$ together with $\left[  P_{n,1}\right]  $ form a set of free
generators of $K_{0}\left(  C\left(  \mathbb{C}P_{q}^{n}\right)  \right)
\cong\mathbb{Z}^{n+1}$.

$\square$

It is of interest to point out that these projections $P_{j,1}$ freely
generating $K_{0}\left(  C\left(  \mathbb{C}P_{q}^{n}\right)  \right)  $ are
actually lying inside $C\left(  \mathbb{C}P_{q}^{n}\right)  \equiv
M_{1}\left(  C\left(  \mathbb{C}P_{q}^{n}\right)  \right)  \subset M_{\infty
}\left(  C\left(  \mathbb{C}P_{q}^{n}\right)  \right)  $ and they form an
increasing finite sequence of projections.

\section{Quantum line bundles over $C\left(  \mathbb{C}P_{q}^{n}\right)  $}

In this section, we identify the quantum line bundles $L_{k}\equiv C\left(
\mathbb{S}_{q}^{2n+1}\right)  _{k}$ of degree $k$ over $C\left(
\mathbb{C}P_{q}^{n}\right)  $ with a concrete (equivalence class of)
projection described in terms of the basic projections. We remark that an
intriguing noncommutative geometric study of these line bundles in comparison
with Adam's classical results on $\mathbb{C}P^{n}$ has been successfully
accomplished by Arici, Brain, and Landi in \cite{AriBrLa}. (The degree
convention is different in the $\pm$-sign.)

To distinguish between ordinary function product and convolution product, we
denote the groupoid C*-algebraic (convolution) multiplication of elements in
$C_{c}\left(  \mathcal{G}\right)  \subset C^{\ast}\left(  \mathcal{G}\right)
$ by $\ast$, while omitting $\ast$ when the elements are represented as
operators or when they are multiplied together pointwise as functions on
$\mathcal{G}$. We also view $C_{c}\left(  \mathfrak{F}_{n}\right)  $ or
$C_{c}\left(  \left(  \mathfrak{F}_{n}\right)  _{k}\right)  $ (also
abbreviated as $C_{c}\left(  \mathfrak{F}_{n}\right)  _{k}$) as left
$C_{c}\left(  \mathfrak{F}_{n}\right)  _{0}$-modules with $C_{c}\left(
\mathfrak{F}_{n}\right)  $ carrying the convolution algebra structure as a
subalgebra of the groupoid C*-algebra $C^{\ast}\left(  \mathfrak{F}%
_{n}\right)  $. Similarly, for a closed subset $X$ of the unit space of
$\mathfrak{F}_{n}$, the inverse image $\mathfrak{F}_{n}\upharpoonright_{X}$ of
$X$ under the source map of $\mathfrak{F}_{n}$ or its grade-$k$ component
$\left(  \mathfrak{F}_{n}\upharpoonright_{X}\right)  _{k}\equiv\left(
\mathfrak{F}_{n}\right)  _{k}\upharpoonright_{X}$ gives rise to a left
$C_{c}\left(  \mathfrak{F}_{n}\right)  _{0}$-module $C_{c}\left(
\mathfrak{F}_{n}\upharpoonright_{X}\right)  $ or $C_{c}\left(  \mathfrak{F}%
_{n}\upharpoonright_{X}\right)  _{k}$.

For $k\in\mathbb{Z}_{\geq}$, the characteristic function $\chi_{B_{k}}\in
C_{c}\left(  \left(  \mathfrak{F}_{n}\right)  _{0}\right)  $ of the compact
open set%
\[
B_{k}:=\left\{  \left[  \left(  0,0^{\left(  n\right)  },w\right)  \right]
\in\left(  \mathfrak{F}_{n}\right)  _{0}:w_{1}\geq k\right\}
\]
is a projection over $C^{\ast}\left(  \left(  \mathfrak{F}_{n}\right)
_{0}\right)  \equiv C\left(  \mathbb{C}P_{q}^{n}\right)  $ which is
represented under $\pi_{n}$ as $P_{-k}\otimes\left(  \otimes^{n-1}I\right)  $,
and
\[
C_{c}\left(  \left(  \mathfrak{F}_{n}\right)  _{0}\right)  \ast\chi_{B_{k}%
}=C_{c}\left(  \left(  \mathfrak{F}_{n}\right)  _{0}\upharpoonright_{B_{k}%
}\right)
\]
where $B_{k}\subset\left(  \mathfrak{F}_{n}\right)  _{0}$ in the notation
$\upharpoonright_{B_{k}}$ is canonically viewed as a subset of the unit space
of $\left(  \mathfrak{F}_{n}\right)  _{0}$.

For $k\leq0$, it is straightforward to check that
\[
\left[  \left(  k,x,w\right)  \right]  \in\left(  \mathfrak{F}_{n}\right)
_{k}\mapsto\left[  \left(  0,x_{1}+k,x_{2},..,x_{n},w_{1}-k,w_{2}%
,..,w_{n}\right)  \right]  \in\left(  \mathfrak{F}_{n}\right)  _{0}%
\upharpoonright_{B_{\left\vert k\right\vert }}%
\]
well defines a bijective homeomorphism. For example, for $w_{1}=\infty$, we
have $x_{1}=-k$ on the domain side and $x_{1}+k=0$ on the range side of this
map, matching the implicit constraints imposed on $\left(  \mathfrak{F}%
_{n}\right)  _{k}$ and $\left(  \mathfrak{F}_{n}\right)  _{0}$. Furthermore
since any $\left[  \left(  k,x,w\right)  \right]  \in\left(  \mathfrak{F}%
_{n}\right)  _{k}$ and its image $\left[  \left(  0,x_{1}+k,x_{2}%
,..,x_{n},w_{1}-k,w_{2},..,w_{n}\right)  \right]  $ share the same target
element $\left[  x+w\right]  \in\overline{\mathbb{Z}}_{\geq}^{n}/\sim$, it
induces a left $C_{c}\left(  \left(  \mathfrak{F}_{n}\right)  _{0}\right)
$-module isomorphism
\[
C_{c}\left(  \left(  \mathfrak{F}_{n}\right)  _{k}\right)  \rightarrow
C_{c}\left(  \left(  \mathfrak{F}_{n}\right)  _{0}\upharpoonright
_{B_{\left\vert k\right\vert }}\right)  \equiv C_{c}\left(  \left(
\mathfrak{F}_{n}\right)  _{0}\right)  \ast\chi_{B_{\left\vert k\right\vert }}%
\]
which extends to a left $C\left(  \mathbb{C}P_{q}^{n}\right)  $-module
isomorphism
\[
L_{k}\equiv\overline{C_{c}\left(  \left(  \mathfrak{F}_{n}\right)
_{k}\right)  }\cong C\left(  \mathbb{C}P_{q}^{n}\right)  \left(
P_{-\left\vert k\right\vert }\otimes\left(  \otimes^{n-1}I\right)  \right)  ,
\]
i.e. the quantum line bundle $L_{k}$ for $k\leq0$ is the finitely generated
left projective module determined by the projection $P_{-\left\vert
k\right\vert }\otimes\left(  \otimes^{n-1}I\right)  $ over $C\left(
\mathbb{C}P_{q}^{n}\right)  $.

For $k>0$, the situation is much more complicated. We first define the closed
open set
\[
\left(  \mathfrak{F}_{n}\right)  _{k,j}:=\left(  \mathfrak{F}_{n}%
\upharpoonright_{\left(  \left\{  0\right\}  ^{j}\times\overline{\mathbb{Z}%
}_{\geq}^{n-j}\right)  /\sim}\right)  _{k}\equiv\left\{  \left[  \left(
k,x,w\right)  \right]  \in\left(  \mathfrak{F}_{n}\right)  _{k}:w\in\left\{
0\right\}  ^{j}\times\overline{\mathbb{Z}}_{\geq}^{n-j}\right\}
\]
with each $C_{c}\left(  \left(  \mathfrak{F}_{n}\right)  _{k,j}\right)  $ a
left $C_{c}\left(  \left(  \mathfrak{F}_{n}\right)  _{0}\right)  $-module.
Note that
\[
C_{c}\left(  \left(  \mathfrak{F}_{n}\right)  _{0,j}\right)  =C_{c}\left(
\left(  \mathfrak{F}_{n}\right)  _{0}\right)  \ast\chi_{\left(  \left\{
0\right\}  ^{j}\times\overline{\mathbb{Z}}_{\geq}^{n-j}\right)  /\sim}%
\]
with $\chi_{\left(  \left\{  0\right\}  ^{j}\times\overline{\mathbb{Z}}_{\geq
}^{n-j}\right)  /\sim}$ represented under $\pi_{n}$ as the projection $\left(
\otimes^{j}P_{1}\right)  \otimes\left(  \otimes^{n-j}I\right)  $ over
$C\left(  \mathbb{C}P_{q}^{n}\right)  $.

Now the left $C_{c}\left(  \left(  \mathfrak{F}_{n}\right)  _{0}\right)
$-module $C_{c}\left(  \left(  \mathfrak{F}_{n}\right)  _{k,j}\right)  $ can
be decomposed as
\[
C_{c}\left(  \left(  \mathfrak{F}_{n}\right)  _{k,j}\right)  =C_{c}\left(
\left(  \mathfrak{F}_{n}\upharpoonright_{\left(  \left\{  0\right\}
^{j}\times\overline{\mathbb{Z}}_{\geq k}\times\overline{\mathbb{Z}}_{\geq
}^{n-j-1}\right)  /\sim}\right)  _{k}\right)  \oplus%
%TCIMACRO{\dbigoplus \limits_{l=0}^{k-1}}%
%BeginExpansion
{\displaystyle\bigoplus\limits_{l=0}^{k-1}}
%EndExpansion
C_{c}\left(  \left(  \mathfrak{F}_{n}\upharpoonright_{\left(  \left\{
0\right\}  ^{j}\times\left\{  l\right\}  \times\overline{\mathbb{Z}}_{\geq
}^{n-j-1}\right)  /\sim}\right)  _{k}\right)  .
\]

It is straightforward to check that
\[
\left[  \left(  k,x,w\right)  \right]  \in\left(  \mathfrak{F}_{n}%
\upharpoonright_{\left(  \left\{  0\right\}  ^{j}\times\overline{\mathbb{Z}%
}_{\geq k}\times\overline{\mathbb{Z}}_{\geq}^{n-j-1}\right)  /\sim}\right)
_{k}\mapsto
\]%
\[
\left[  \left(  0,x_{1},..,x_{j},x_{j+1}+k,x_{j+2},..,x_{n},0^{\left(
j\right)  },w_{j+1}-k,w_{j+2},..,w_{n}\right)  \right]  \in\left(
\mathfrak{F}_{n}\upharpoonright_{\left(  \left\{  0\right\}  ^{j}%
\times\overline{\mathbb{Z}}_{\geq}^{n-j}\right)  /\sim}\right)  _{0}%
\equiv\left(  \mathfrak{F}_{n}\right)  _{0,j}%
\]
well defines a bijective homeomorphism. For example, we are considering only
$w$ with $w_{1}=...=w_{j}=0<\infty$, while for $w_{j+1}=\infty$, we have
$x_{j+1}=-k-x_{1}-\cdots-x_{j}$ on the domain side and $-x_{1}-\cdots
-x_{j}-\left(  x_{j+1}+k\right)  =0$ on the range side of this map, matching
the implicit constraints imposed on $\left(  \mathfrak{F}_{n}\right)  _{k}$
and $\left(  \mathfrak{F}_{n}\right)  _{0}$. Furthermore since any $\left[
\left(  k,x,w\right)  \right]  \ $ and its image under this bijection share
the same target element $\left[  x+w\right]  \in\overline{\mathbb{Z}}_{\geq
}^{n}/\sim$, it induces a left $C_{c}\left(  \left(  \mathfrak{F}_{n}\right)
_{0}\right)  $-module isomorphism
\[
C_{c}\left(  \left(  \mathfrak{F}_{n}\upharpoonright_{\left(  \left\{
0\right\}  ^{j}\times\overline{\mathbb{Z}}_{\geq k}\times\overline{\mathbb{Z}%
}_{\geq}^{n-j-1}\right)  /\sim}\right)  _{k}\right)  \rightarrow C_{c}\left(
\left(  \mathfrak{F}_{n}\right)  _{0,j}\right)  =C_{c}\left(  \left(
\mathfrak{F}_{n}\right)  _{0}\right)  \ast\chi_{\left(  \left\{  0\right\}
^{j}\times\overline{\mathbb{Z}}_{\geq}^{n-j}\right)  /\sim}%
\]
which extends to a left $C\left(  \mathbb{C}P_{q}^{n}\right)  $-module
isomorphism%
\[
\overline{C_{c}\left(  \left(  \mathfrak{F}_{n}\upharpoonright_{\left(
\left\{  0\right\}  ^{j}\times\overline{\mathbb{Z}}_{\geq k}\times
\overline{\mathbb{Z}}_{\geq}^{n-j-1}\right)  /\sim}\right)  _{k}\right)
}\cong C\left(  \mathbb{C}P_{q}^{n}\right)  \left(  \left(  \otimes^{j}%
P_{1}\right)  \otimes\left(  \otimes^{n-j}I\right)  \right)  .
\]

On the other hand, for any $0\leq l\leq k-1$,
\[
\left[  \left(  k,x,w\right)  \right]  \in\left(  \mathfrak{F}_{n}%
\upharpoonright_{\left(  \left\{  0\right\}  ^{j}\times\left\{  l\right\}
\times\overline{\mathbb{Z}}_{\geq}^{n-j-1}\right)  /\sim}\right)  _{k}\mapsto
\]%
\[
\left[  \left(  k-l,x_{1},..,x_{j},x_{j+1}+l,x_{j+2},..,x_{n},0^{\left(
j+1\right)  },w_{j+2},..,w_{n}\right)  \right]  \in\left(  \mathfrak{F}%
_{n}\upharpoonright_{\left(  \left\{  0\right\}  ^{j+1}\times\overline
{\mathbb{Z}}_{\geq}^{n-j-1}\right)  /\sim}\right)  _{k-l}\equiv\left(
\mathfrak{F}_{n}\right)  _{k-l,j+1}%
\]
well defines a bijective homeomorphism which preserves the target element
$\left[  x+w\right]  $ and hence induces a left $C_{c}\left(  \left(
\mathfrak{F}_{n}\right)  _{0}\right)  $-module isomorphism
\[
C_{c}\left(  \left(  \mathfrak{F}_{n}\upharpoonright_{\left(  \left\{
0\right\}  ^{j}\times\left\{  l\right\}  \times\overline{\mathbb{Z}}_{\geq
}^{n-j-1}\right)  /\sim}\right)  _{k}\right)  \rightarrow C_{c}\left(  \left(
\mathfrak{F}_{n}\right)  _{k-l,j+1}\right)  .
\]

So summarizing, we get the isomorphism relation%
\[
\text{(*)\ \ \ }C_{c}\left(  \left(  \mathfrak{F}_{n}\right)  _{k,j}\right)
\cong C_{c}\left(  \left(  \mathfrak{F}_{n}\right)  _{0,j}\right)  \oplus%
%TCIMACRO{\dbigoplus \limits_{l=0}^{k-1}}%
%BeginExpansion
{\displaystyle\bigoplus\limits_{l=0}^{k-1}}
%EndExpansion
C_{c}\left(  \left(  \mathfrak{F}_{n}\right)  _{k-l,j+1}\right)  \equiv
C_{c}\left(  \left(  \mathfrak{F}_{n}\right)  _{0,j}\right)  \oplus%
%TCIMACRO{\dbigoplus \limits_{l=1}^{k}}%
%BeginExpansion
{\displaystyle\bigoplus\limits_{l=1}^{k}}
%EndExpansion
C_{c}\left(  \left(  \mathfrak{F}_{n}\right)  _{l,j+1}\right)
\]
which is recursive in the sense that the right hand side contains terms with
either $k$ decreased or $j$ increased. So repeated application of this
recursive expansion can lead to a direct sum of terms of the form
$C_{c}\left(  \left(  \mathfrak{F}_{n}\right)  _{0,m}\right)  $ or the form
$C_{c}\left(  \left(  \mathfrak{F}_{n}\right)  _{l,n}\right)  $, where
\[
\overline{C_{c}\left(  \left(  \mathfrak{F}_{n}\right)  _{0,m}\right)  }\cong
C\left(  \mathbb{C}P_{q}^{n}\right)  \left(  \left(  \otimes^{m}P_{1}\right)
\otimes\left(  \otimes^{n-m}I\right)  \right)
\]
while
\[
\left[  \left(  l,x,w\right)  \right]  \equiv\left[  \left(  l,x,0^{\left(
n\right)  }\right)  \right]  \in\left(  \mathfrak{F}_{n}\upharpoonright
_{\left\{  0\right\}  ^{n}/\sim}\right)  _{l}\equiv\left(  \mathfrak{F}%
_{n}\right)  _{l,n}\mapsto\left[  \left(  0,x,0^{\left(  n\right)  }\right)
\right]  \in\left(  \mathfrak{F}_{n}\right)  _{0,n}%
\]
well defines a bijective homeomorphism which induces a left $C_{c}\left(
\left(  \mathfrak{F}_{n}\right)  _{0}\right)  $-module isomorphism
\[
C_{c}\left(  \left(  \mathfrak{F}_{n}\right)  _{l,n}\right)  \rightarrow
C_{c}\left(  \left(  \mathfrak{F}_{n}\right)  _{0,n}\right)
\]
extending to a left $C\left(  \mathbb{C}P_{q}^{n}\right)  $-module
isomorphism
\[
\overline{C_{c}\left(  \left(  \mathfrak{F}_{n}\right)  _{l,n}\right)  }\cong
C\left(  \mathbb{C}P_{q}^{n}\right)  \left(  \otimes^{n}P_{1}\right)  .
\]

\textbf{Theorem 3}. For $n\geq1$, the quantum line bundle $L_{k}\equiv
C\left(  \mathbb{S}_{q}^{2n+1}\right)  _{k}$ of degree $k\in\mathbb{Z}$ over
$C\left(  \mathbb{C}P_{q}^{n}\right)  $ is isomorphic to the finitely
generated projective left module over $C\left(  \mathbb{C}P_{q}^{n}\right)  $
determined by the projection $P_{-\left\vert k\right\vert }\otimes\left(
\otimes^{n-1}I\right)  $ if $k\leq0$ (with $P_{-0}:=I$ understood), and the
projection%
\[
\boxplus_{j=0}^{n}\left(  \boxplus^{C_{j}^{k+j-1}}\left(  \left(  \otimes
^{j}P_{1}\right)  \otimes\left(  \otimes^{n-j}I\right)  \right)  \right)
\]
if $k>0$, where $C_{j}^{k}$ denotes the combinatorial number $\left(
k!\right)  /\left(  j!\left(  k-j\right)  !\right)  $.

Proof. Having already taken care of the case of $k\leq0$ in the above
discussion, we only need to consider the case of $k>0$.

First we establish by induction on $l$ that
\[
\text{(**)\ \ }C_{c}\left(  \left(  \mathfrak{F}_{n}\right)  _{k,0}\right)
\cong\left(
%TCIMACRO{\dbigoplus \limits_{j=0}^{l-1}}%
%BeginExpansion
{\displaystyle\bigoplus\limits_{j=0}^{l-1}}
%EndExpansion
\left(  \oplus^{C_{j}^{k+j-1}}C_{c}\left(  \left(  \mathfrak{F}_{n}\right)
_{0,j}\right)  \right)  \right)  \oplus\left(
%TCIMACRO{\dbigoplus \limits_{m=1}^{k}}%
%BeginExpansion
{\displaystyle\bigoplus\limits_{m=1}^{k}}
%EndExpansion
\left(  \oplus^{C_{l-1}^{k-m+l-1}}C_{c}\left(  \left(  \mathfrak{F}%
_{n}\right)  _{m,l}\right)  \right)  \right)  .
\]
Indeed for $l=1$, (**) becomes
\[
C_{c}\left(  \left(  \mathfrak{F}_{n}\right)  _{k,0}\right)  \cong
C_{c}\left(  \left(  \mathfrak{F}_{n}\right)  _{0,0}\right)  \oplus\left(
%TCIMACRO{\dbigoplus \limits_{m=1}^{k}}%
%BeginExpansion
{\displaystyle\bigoplus\limits_{m=1}^{k}}
%EndExpansion
C_{c}\left(  \left(  \mathfrak{F}_{n}\right)  _{m,1}\right)  \right)  ,
\]
which is the same as the established recursive relation (*) with $j=0$. For
$n\geq l>1$, by the induction hypothesis for $l-1$ and the recursive relation
(*), we get
\[
C_{c}\left(  \left(  \mathfrak{F}_{n}\right)  _{k,0}\right)  \cong\left(
%TCIMACRO{\dbigoplus \limits_{j=0}^{l-2}}%
%BeginExpansion
{\displaystyle\bigoplus\limits_{j=0}^{l-2}}
%EndExpansion
\left(  \oplus^{C_{j}^{k+j-1}}C_{c}\left(  \left(  \mathfrak{F}_{n}\right)
_{0,j}\right)  \right)  \right)  \oplus\left(
%TCIMACRO{\dbigoplus \limits_{m=1}^{k}}%
%BeginExpansion
{\displaystyle\bigoplus\limits_{m=1}^{k}}
%EndExpansion
\left(  \oplus^{C_{l-2}^{k-m+l-2}}C_{c}\left(  \left(  \mathfrak{F}%
_{n}\right)  _{m,l-1}\right)  \right)  \right)
\]%
\[
\cong\left(
%TCIMACRO{\dbigoplus \limits_{j=0}^{l-2}}%
%BeginExpansion
{\displaystyle\bigoplus\limits_{j=0}^{l-2}}
%EndExpansion
\left(  \oplus^{C_{j}^{k+j-1}}C_{c}\left(  \left(  \mathfrak{F}_{n}\right)
_{0,j}\right)  \right)  \right)  \oplus\left(
%TCIMACRO{\dbigoplus \limits_{m=1}^{k}}%
%BeginExpansion
{\displaystyle\bigoplus\limits_{m=1}^{k}}
%EndExpansion
\left(  \oplus^{C_{l-2}^{k-m+l-2}}\left(  C_{c}\left(  \left(  \mathfrak{F}%
_{n}\right)  _{0,l-1}\right)  \oplus%
%TCIMACRO{\dbigoplus \limits_{i=1}^{m}}%
%BeginExpansion
{\displaystyle\bigoplus\limits_{i=1}^{m}}
%EndExpansion
C_{c}\left(  \left(  \mathfrak{F}_{n}\right)  _{i,l}\right)  \right)  \right)
\right)
\]%
\[
\cong\left(
%TCIMACRO{\dbigoplus \limits_{j=0}^{l-2}}%
%BeginExpansion
{\displaystyle\bigoplus\limits_{j=0}^{l-2}}
%EndExpansion
\left(  \oplus^{C_{j}^{k+j-1}}C_{c}\left(  \left(  \mathfrak{F}_{n}\right)
_{0,j}\right)  \right)  \right)  \oplus\left(  \oplus^{%
%TCIMACRO{\dsum \limits_{m=1}^{k}}%
%BeginExpansion
{\displaystyle\sum\limits_{m=1}^{k}}
%EndExpansion
C_{l-2}^{k-m+l-2}}C_{c}\left(  \left(  \mathfrak{F}_{n}\right)  _{0,l-1}%
\right)  \right)  \oplus\left(
%TCIMACRO{\dbigoplus \limits_{m=1}^{k}}%
%BeginExpansion
{\displaystyle\bigoplus\limits_{m=1}^{k}}
%EndExpansion%
%TCIMACRO{\dbigoplus \limits_{i=1}^{m}}%
%BeginExpansion
{\displaystyle\bigoplus\limits_{i=1}^{m}}
%EndExpansion
\left(  \oplus^{C_{l-2}^{k-m+l-2}}C_{c}\left(  \left(  \mathfrak{F}%
_{n}\right)  _{i,l}\right)  \right)  \right)
\]%
\[
\cong\left(
%TCIMACRO{\dbigoplus \limits_{j=0}^{l-2}}%
%BeginExpansion
{\displaystyle\bigoplus\limits_{j=0}^{l-2}}
%EndExpansion
\left(  \oplus^{C_{j}^{k+j-1}}C_{c}\left(  \left(  \mathfrak{F}_{n}\right)
_{0,j}\right)  \right)  \right)  \oplus\left(  \oplus^{C_{l-1}^{k+l-2}}%
C_{c}\left(  \left(  \mathfrak{F}_{n}\right)  _{0,l-1}\right)  \right)
\oplus\left(
%TCIMACRO{\dbigoplus \limits_{i=1}^{k}}%
%BeginExpansion
{\displaystyle\bigoplus\limits_{i=1}^{k}}
%EndExpansion%
%TCIMACRO{\dbigoplus \limits_{m=i}^{k}}%
%BeginExpansion
{\displaystyle\bigoplus\limits_{m=i}^{k}}
%EndExpansion
\left(  \oplus^{C_{l-2}^{k-m+l-2}}C_{c}\left(  \left(  \mathfrak{F}%
_{n}\right)  _{i,l}\right)  \right)  \right)
\]%
\[
\cong\left(
%TCIMACRO{\dbigoplus \limits_{j=0}^{l-1}}%
%BeginExpansion
{\displaystyle\bigoplus\limits_{j=0}^{l-1}}
%EndExpansion
\left(  \oplus^{C_{j}^{k+j-1}}C_{c}\left(  \left(  \mathfrak{F}_{n}\right)
_{0,j}\right)  \right)  \right)  \oplus\left(
%TCIMACRO{\dbigoplus \limits_{i=1}^{k}}%
%BeginExpansion
{\displaystyle\bigoplus\limits_{i=1}^{k}}
%EndExpansion
\left(  \oplus^{\sum_{m=i}^{k}C_{l-2}^{k-m+l-2}}C_{c}\left(  \left(
\mathfrak{F}_{n}\right)  _{i,l}\right)  \right)  \right)
\]%
\[
\cong\left(
%TCIMACRO{\dbigoplus \limits_{j=0}^{l-1}}%
%BeginExpansion
{\displaystyle\bigoplus\limits_{j=0}^{l-1}}
%EndExpansion
\left(  \oplus^{C_{j}^{k+j-1}}C_{c}\left(  \left(  \mathfrak{F}_{n}\right)
_{0,j}\right)  \right)  \right)  \oplus\left(
%TCIMACRO{\dbigoplus \limits_{i=1}^{k}}%
%BeginExpansion
{\displaystyle\bigoplus\limits_{i=1}^{k}}
%EndExpansion
\left(  \oplus^{C_{l-1}^{k-i+l-1}}C_{c}\left(  \left(  \mathfrak{F}%
_{n}\right)  _{i,l}\right)  \right)  \right)
\]%
\[
\equiv\left(
%TCIMACRO{\dbigoplus \limits_{j=0}^{l-1}}%
%BeginExpansion
{\displaystyle\bigoplus\limits_{j=0}^{l-1}}
%EndExpansion
\left(  \oplus^{C_{j}^{k+j-1}}C_{c}\left(  \left(  \mathfrak{F}_{n}\right)
_{0,j}\right)  \right)  \right)  \oplus\left(
%TCIMACRO{\dbigoplus \limits_{m=1}^{k}}%
%BeginExpansion
{\displaystyle\bigoplus\limits_{m=1}^{k}}
%EndExpansion
\left(  \oplus^{C_{l-1}^{k-m+l-1}}C_{c}\left(  \left(  \mathfrak{F}%
_{n}\right)  _{m,l}\right)  \right)  \right)
\]
where
\[%
%TCIMACRO{\dsum \limits_{m=1}^{k}}%
%BeginExpansion
{\displaystyle\sum\limits_{m=1}^{k}}
%EndExpansion
C_{l-2}^{k-m+l-2}=C_{l-2}^{l-2}+C_{l-2}^{l-1}+C_{l-2}^{l}+C_{l-2}^{l+1}%
+\cdots+C_{l-2}^{k+l-3}%
\]%
\[
=C_{l-1}^{l-1}+C_{l-2}^{l-1}+C_{l-2}^{l}+C_{l-2}^{l+1}+\cdots+C_{l-2}^{k+l-3}%
\]%
\[
=C_{l-1}^{l}+C_{l-2}^{l}+C_{l-2}^{l+1}+\cdots+C_{l-2}^{k+l-3}=C_{l-1}%
^{l+1}+C_{l-2}^{l+1}+\cdots+C_{l-2}^{k+l-3}%
\]%
\[
=\cdots=C_{l-1}^{k+l-3}+C_{l-2}^{k+l-3}=C_{l-1}^{k+l-2}%
\]
and similarly%
\[
\sum_{m=i}^{k}C_{l-2}^{k-m+l-2}=C_{l-2}^{l-2}+C_{l-2}^{l-1}+\cdots
+C_{l-2}^{k-i+l-2}=C_{l-1}^{k-i+l-1}.
\]
Thus (**) holds for $n\geq l>1$, concluding the inductive proof of (**).

Now by (**) for $l=n$, we get
\[
L_{k}\equiv\overline{C_{c}\left(  \left(  \mathfrak{F}_{n}\right)
_{k,0}\right)  }\cong\left(
%TCIMACRO{\dbigoplus \limits_{j=0}^{n-1}}%
%BeginExpansion
{\displaystyle\bigoplus\limits_{j=0}^{n-1}}
%EndExpansion
\left(  \oplus^{C_{j}^{k+j-1}}\overline{C_{c}\left(  \left(  \mathfrak{F}%
_{n}\right)  _{0,j}\right)  }\right)  \right)  \oplus\left(
%TCIMACRO{\dbigoplus \limits_{m=1}^{k}}%
%BeginExpansion
{\displaystyle\bigoplus\limits_{m=1}^{k}}
%EndExpansion
\left(  \oplus^{C_{n-1}^{k-m+n-1}}\overline{C_{c}\left(  \left(
\mathfrak{F}_{n}\right)  _{m,n}\right)  }\right)  \right)
\]%
\[
=\left(
%TCIMACRO{\dbigoplus \limits_{j=0}^{n-1}}%
%BeginExpansion
{\displaystyle\bigoplus\limits_{j=0}^{n-1}}
%EndExpansion
\left(  \oplus^{C_{j}^{k+j-1}}C\left(  \mathbb{C}P_{q}^{n}\right)  \left(
\left(  \otimes^{j}P_{1}\right)  \otimes\left(  \otimes^{n-j}I\right)
\right)  \right)  \right)  \oplus\left(
%TCIMACRO{\dbigoplus \limits_{m=1}^{k}}%
%BeginExpansion
{\displaystyle\bigoplus\limits_{m=1}^{k}}
%EndExpansion
\left(  \oplus^{C_{n-1}^{k-m+n-1}}C\left(  \mathbb{C}P_{q}^{n}\right)  \left(
\otimes^{n}P_{1}\right)  \right)  \right)
\]%
\[
=\left(
%TCIMACRO{\dbigoplus \limits_{j=0}^{n-1}}%
%BeginExpansion
{\displaystyle\bigoplus\limits_{j=0}^{n-1}}
%EndExpansion
\left(  \oplus^{C_{j}^{k+j-1}}C\left(  \mathbb{C}P_{q}^{n}\right)  \left(
\left(  \otimes^{j}P_{1}\right)  \otimes\left(  \otimes^{n-j}I\right)
\right)  \right)  \right)  \oplus\left(  \oplus^{\sum_{m=1}^{k}C_{n-1}%
^{k-m+n-1}}C\left(  \mathbb{C}P_{q}^{n}\right)  \left(  \otimes^{n}%
P_{1}\right)  \right)
\]%
\[
=\left(
%TCIMACRO{\dbigoplus \limits_{j=0}^{n-1}}%
%BeginExpansion
{\displaystyle\bigoplus\limits_{j=0}^{n-1}}
%EndExpansion
\left(  \oplus^{C_{j}^{k+j-1}}C\left(  \mathbb{C}P_{q}^{n}\right)  \left(
\left(  \otimes^{j}P_{1}\right)  \otimes\left(  \otimes^{n-j}I\right)
\right)  \right)  \right)  \oplus\left(  \oplus^{C_{n}^{k+n-1}}C\left(
\mathbb{C}P_{q}^{n}\right)  \left(  \otimes^{n}P_{1}\right)  \right)
\]%
\[
=%
%TCIMACRO{\dbigoplus \limits_{j=0}^{n}}%
%BeginExpansion
{\displaystyle\bigoplus\limits_{j=0}^{n}}
%EndExpansion
\left(  \oplus^{C_{j}^{k+j-1}}C\left(  \mathbb{C}P_{q}^{n}\right)  \left(
\left(  \otimes^{j}P_{1}\right)  \otimes\left(  \otimes^{n-j}I\right)
\right)  \right)
\]
where again%
\[
\sum_{m=1}^{k}C_{n-1}^{k-m+n-1}=C_{n-1}^{n-1}+C_{n-1}^{n}+C_{n-1}^{n+1}%
+\cdots+C_{n-1}^{k+n-2}=C_{n}^{k+n-1}.
\]

Thus $L_{k}$ for $k>0$ is implemented by the projection
\[
\boxplus_{j=0}^{n}\left(  \boxplus^{C_{j}^{k+j-1}}\left(  \left(  \otimes
^{j}P_{1}\right)  \otimes\left(  \otimes^{n-j}I\right)  \right)  \right)  .
\]

$\square$

\end{document}